\documentclass[12pt]{article}
\usepackage[english]{babel}
\usepackage{epsfig}
\usepackage{amsfonts}
\usepackage{amssymb}
\usepackage{amsmath}
\usepackage{amsthm}
\usepackage{latexsym}
\usepackage{graphicx}
\usepackage{bm}
\usepackage{tabularx}
\usepackage{booktabs} 
\usepackage{enumerate}
\usepackage{caption}
\usepackage{wrapfig}  

\usepackage{mathtools}  
\usepackage{eqnarray}
\usepackage{enumitem}

\usepackage[titletoc]{appendix}

\usepackage[numbers]{natbib}
\usepackage{comment}

\usepackage{bbm} 
\usepackage{url}
\usepackage{subfig}
\usepackage{hyperref}
\usepackage{color}
\usepackage{float}

\usepackage{hyperref}
\hypersetup{colorlinks=true, citecolor=blue, linkcolor=red, urlcolor=green}
 
\usepackage{bbm}
\usepackage{url}
\usepackage{color}
\usepackage{xcolor}
\usepackage[section]{placeins}

\theoremstyle{plain}
\newtheorem{theorem}{Theorem}[section]
\newtheorem{lemma}[theorem]{Lemma}
\newtheorem{corollary}[theorem]{Corollary}
\newtheorem{proposition}[theorem]{Proposition}

\theoremstyle{definition}

\theoremstyle{remark}

\newcommand{\shortdot}[1]{\raisebox{-0.4pt}{$\stackrel{\bullet}{#1}$}}

\title{Overlap Times in the Infinite Server Queue}
\author{ Sergio Palomo   \\Systems Engineering \\ Cornell
University \\ {sdp85@cornell.edu} 
\and
Jamol Pender \footnote{Corresponding Author}
\\
School of Operations Research and Information Engineering
\\
Cornell University
\\
{jjp274@cornell.edu}
}

\usepackage{natbib}
\usepackage{graphicx}
\usepackage{amsmath}

\begin{document}

\maketitle

\begin{abstract}
Imagine, you enter a grocery store to buy food.  How many people do you overlap with in this store? How much time do you overlap with each person in the store?  In this paper, we answer these questions by studying the overlap times between customers in the infinite server queue. We compute in closed form the steady state distribution of the overlap time between a pair of customers and the distribution of the number of customers that an arriving customer will overlap with.  Finally, we define a residual process that counts the number of overlapping customers that overlap in the queue for at least $\delta$ time units and compute its mean, variance, and distribution in the exponential service setting.  
\end{abstract}

\section{Introduction}
Since the beginning of the COVID-19 pandemic, much work has focused on using deterministic compartmentalized style models to estimate the infection rate and dynamics of the spread, see for example \citet{dandekar2021safe, nguemdjo2020simulating, kaplan2020om}.  However, we know that stochastic effects can play an important role in determining the spread, see for example \citet{drakopoulos2017network, palomo2020flattening, pang2020functional, forien2020epidemic, moein2021inefficiency}.  As shoppers crowd a store to stock up on water or large amounts of non-perishable items it is inevitable that the virus would spread.  To combat the spread of the virus, many service facilities and systems have installed new air filters, transparent barriers and have asked that patrons wear facial masks.  Moreover, these service systems also have implemented various forms of social/physical distancing in order to minimize close proximity of one customer to another \citet{bove2020restrict}.  

However, there are some places where people work or shop that limiting distance is not feasible.  In this case, we really care about how much customers overlap with one another.   Recently there has been new work by \citet{kang2021queueing} and \citet{palomooverlap} that explores how one can calculate the overlap times of customers in single server queues.  More specifically, \citet{kang2021queueing} shows how to use the overlaps to compute a new $R_0$ value for understanding infection rates in compartmentalized epidemic models and \citet{palomooverlap} proves that the overlap distribution is exponential for the M/M/1 and shows via simulation that a similar result holds for the non-Markovian setting as well.  Our analysis is important because it can demonstrate, exactly, how much overlap occurs and can provide distributional information or prediction intervals for possible overlap.  Moreover, it can be used as a tool to prevent large overlaps and our analysis can be used as a design tool construct appropriate overlap by restricting the arrival rate or service distributions. 

In this paper, we extend the overlap time analysis to the setting of infinite server queue.  At first glance the infinite server queue analysis might not seem relevant, however, for service entities such as grocery stores, outlets, restaurants, and retail shops, an infinite server queue is quite relevant as there is not much waiting or the waiting to check out may be insignificant when compared to the shopping time.  Moreover, the overlap times in the infinite server queue serve as a lower bound for the overlap times that a customer might experience in systems where there is significant waiting or in a multi-server setting.  What also makes the infinite server queue important is that we are able to derive explicit formulas for the overlap distribution and residual overlap distributions as well as the number of people that a customer will overlap with during the duration of their service experience.  In what follows we describe the contributions of our work and how the rest of the paper is organized.  

\subsection{Contributions of the Paper}
\begin{itemize}
    \item We derive the steady state distribution for the overlap time for customers that are exactly $k$ spaces apart.  
    \item We derive the distribution of the number of customers that any customer overlaps with during their service.  
    \item We also construct a residual overlap process and compute its mean, variance and distribution.  
    \item We use simulation to verify our results.  
\end{itemize}

\subsection{Organization of Paper}
In Section \ref{Stoc:Model}, we describe the stochastic model that we will use in this work.  We derive an equation for describing the overlap times for customers in the infinite server queue.  We use this equation to compute the steady state distribution of the overlap time of customers that are $k$ spaces apart.  In Section \ref{Sec:Overlaps}, we compute the mean and variance of the number of people a customer will overlap with during their time in the queue. We also compute a residual version where a customer must overlap at least $\delta$ units of time and compute its mean and variance as well.  Finally in Section \ref{conclusion}, we provide a conclusion and some future research directions.

\section{Infinite Server Overlap Times} \label{Stoc:Model}

In this section, we study the infinite server queue with the intention of understanding how much time do adjacent customers spend in the system together.  A similar type of analysis has been completed \citet{kang2021queueing, palomooverlap}.  Before we get into studying the overlap times, we review some concepts and results for the $M_t/G/\infty$ queueing model as this will be very helpful in the analysis that follows.  

\subsection{Review of the $M_t/G/\infty$ Queue}

In this section, we review the $M_t/G/\infty$ queueing model. The $M_t/G/\infty$ queue $Q^{\infty}(t)$ has a Poisson distribution with time varying mean $q^{\infty}(t)$.  As observed in \citet{eick1993mt, eick1993physics}, $q^{\infty}(t)$ has the following integral representation
\begin{eqnarray}
q^{\infty}(t) &=& E[Q^{\infty}(t) ] \\
&=& \int^{t}_{-\infty} \overline{G}(t-u) \lambda(u) du \\
&=& E\left[\int^{t}_{t-S} \lambda(u) du \right] \\
&=& E[ \lambda( t - S_e ) ] \cdot E[S]
\end{eqnarray}
where $\lambda(u)$ is the time varying arrival rate and $S$ represents a service time with distribution G, $\overline{G} = 1 - G(t) = 
\mathbb{P}( S > t)$, and $S_e$ is a random variable with distribution that follows the 
stationary excess of residual-lifetime cdf $G_e$, defined by
 \begin{eqnarray}
G_e(t) &\equiv&  \mathbb{P}( S_e < t) = \frac{1}{E[S]} \int^{t}_{0} \overline{G}(u) du = \frac{1}{E[S]} \int^{t}_{0} \mathbb{P}( S > u) du, \ 
\ \ t \geq 0.
\end{eqnarray}

We find it also useful to compute the mean of the $M/M/\infty$ queue using the solution of a linear differential equation.  

\begin{lemma} \label{ode_soln}
Let $q(t)$ be the solution to the following differential equation
\begin{equation}
\shortdot{q} = \lambda - \mu q(t) 
\end{equation}
where $q(0) = q_0$.  Then the solution for any value of $t$ is given by 
\begin{equation}
q(t) = q_0 e^{-\mu t} + \frac{\lambda}{\mu} \left( 1 - e^{-\mu t} \right).
\end{equation}
\begin{proof}
This follows from standard results on ordinary differential equations.  
\end{proof}
\end{lemma}

\subsection{Review of the Incomplete Gamma Function}

Before we begin, we would like to define a few functions that will be useful in our future analysis.  First we will recall the integral definition of the gamma function i.e.
\begin{eqnarray*}
\Gamma(a)  &=&  \int^{\infty}_{0} t^{a-1} e^{-t} dt \quad \mathrm{for} \ a >  0.
\end{eqnarray*}
One can also partition the gamma function integral at a point $x\geq0$ to obtain the lower and upper incomplete gamma functions respectively as
\begin{eqnarray*}
\gamma(a,x)  &=&  \int^{x}_{0} t^{a-1} e^{-t} dt,
\end{eqnarray*}
\begin{eqnarray*}
\Gamma(a,x)  &=&  \int^{\infty}_{x} t^{a-1} e^{-t} dt =  e^{-x} \int^{\infty}_{x} (t + x)^{a-1} e^{-t} dt .
\end{eqnarray*}

Moreover, we know that 
\begin{eqnarray*}
\frac{1}{c^a} \gamma(a,cx)  &=&  \int^{x}_{0} t^{a-1} e^{-ct} dt,
\end{eqnarray*}
\begin{eqnarray*}
\frac{1}{c^a} \Gamma(a,cx)  &=&  \int^{\infty}_{x} t^{a-1} e^{-ct} dt .
\end{eqnarray*}

The lower and upper incomplete gamma functions will be used extensively in what follows.  

\subsection{Overlap Time Distribution}

In this section, we consider the GI/GI/$\infty$ queue.  In this queue we let $A_i$ be the arrival time of the $i^{th}$ customer and we define the inter-arrival time between the $i^{th}$ and $(i+1)^{th}$ customers to be $A_{i+1} - A_i$, which are i.i.d random variables with cumulative distribution function (cdf) $F(x)$.  We also assume that $S_i$ is the service time of the $i^{th}$ customer and the service times are i.i.d with cdf $G(x)$.   In the infinite server queue by definition, no customer will wait.  Thus, the departure time for the $n^{th}$ customer is given by the following equation
\begin{eqnarray}
D_n &=& S_n + A_{n}.
\end{eqnarray}
Now that we understand the arrival and departure time for each customer, we can now construct an equation for the overlap time between adjacent customers.  The overlap time between the $n^{th}$ and $(n+k)^{th}$ customers is given by
\begin{eqnarray}
O_{n,n+k} &=& \left( \min( D_n , D_{n+k} ) - A_{n+k} \right)^+ \\
&=& \left( \min( A_n + S_n , A_{n+k} + S_{n+k} ) - A_{n+k} \right)^+ \\
&=&  \left( (D_n - A_{n+k} )^+ \wedge S_{n+k} \right) \\
&=&  \left( (S_n + A_{n} - A_{n+k} )^+ \wedge S_{n+k} \right) \\
&=&  \left(  (S_n - \left( A_{n+k} - A_{n} \right) )^+ \wedge S_{n+k} \right) .
\end{eqnarray}

It is important to observe that the overlap time between the $n^{th}$ and $(n+k)^{th}$ customers can be decomposed into two parts.  The first part (the left term in the minimum) is the time that the $n^{th}$ customer overlaps with the $(n+k)^{th}$ customer given that the $n^{th}$ customer stays longer.  The second part is the service time of the $(n+k)^{th}$ customer if the service time of the $(n+k)^{th}$ stay is shorter than the $n^{th}$ customers service time minus the inter-arrival time gap.  We will leverage this representation when considering the steady state overlap time and the fact that all of the random variables are independent from one another.  However, before we derive the steady state distribution of the overlap time for the $M/M/\infty$ queue, we find it useful to derive an important lemma about the distribution of an exponential random variable minus an Erlang random variable.  This lemma is given below in Lemma \ref{diff_Z}.

\begin{lemma} \label{diff_Z}
Let $X$ be an Exp($\mu$) random variable  and $Y$ be an Erlang$(k, \lambda)$ random variable.  Then, the probability density function of $Z = X - Y$ is equal to
\begin{equation} \label{diff_dist}
f_Z(z) =
\begin{cases}
1 - \frac{\mu \lambda^k}{(\lambda + \mu)^k} e^{-\lambda z} \sum^{k-1}_{j=0} \frac{ [(\lambda + \mu) z]^j }{j!}  , \quad \mathrm{for} \ z \leq 0 \\ 
 \frac{\mu \lambda^k}{(\lambda + \mu)^k} e^{-\mu z}   \quad \mathrm{for} \ z > 0 
\end{cases}
\end{equation}

 \begin{proof}
See \citet{palomo2020flattening}.   
 \end{proof}
\end{lemma}

\begin{theorem} \label{ssdist}
Let $O_k$ be the steady state distribution of $O_{n,n+k}$ in the $M/M/\infty$ queue, then the tail distribution of $O_k = \lim_{n \to \infty} O_{n,n+k}$ is given by
\begin{eqnarray} \label{overlap_dist}
\mathbb{P} \left( O_{k} > t  \right) &=&  \left( \frac{\lambda}{\lambda + \mu} \right)^k e^{-2\mu t} 
\end{eqnarray}
and 
\begin{eqnarray}
\mathbb{P} \left( O_{k} =  0  \right) &=& 1 -  \left( \frac{\lambda}{\lambda + \mu} \right)^k  .
\end{eqnarray}
\begin{proof}
First, we need to decompose the overlap probability into two probabilities by using a property of the minimum of two independent random variables i.e.
\begin{eqnarray}
\mathbb{P} \left( O_{k} > t  \right) &=& \mathbb{P} \left( \left(  (\mathcal{S} - \mathcal{A}_k )^+ \wedge \tilde{\mathcal{S}} \right)  > t  \right) \\
&=& \mathbb{P} \left( \left(  (\mathcal{S} - \mathcal{A}_k )^+ \right)  > t  \right) \cdot  \mathbb{P} \left(  \tilde{\mathcal{S}}  > t  \right) \\
&=& \mathbb{P} \left( \left(  (\mathcal{S} - \mathcal{A}_k )^+ \right)  > t  \right) \cdot  e^{-\mu t} .
\end{eqnarray}

Now it suffices to compute the remaining probability expression.  We can do this by exploiting Proposition \ref{diff_Z} as follows.  
\begin{eqnarray}
\mathbb{P} \left( O_{k} > t  \right) &=& \mathbb{P} \left( (\mathcal{S} -  \mathcal{A}_k )^+ > t \right) \cdot  e^{-\mu t}  \\
&=& \mathbb{P} \left( \mathcal{S} -  \mathcal{A}_k  > t \right) \cdot  e^{-\mu t}  \\
&=& \mathbb{P} \left( \mathcal{S} > t +  \mathcal{A}_k   \right) \cdot  e^{-\mu t} \\
&=& \left( \int^{\infty}_{t}  \frac{\mu \lambda^k}{(\lambda + \mu)^k} e^{-\mu z} dz  \right) \cdot  e^{-\mu t} \\
&=& \left( \frac{\mu \lambda^k}{(\lambda + \mu)^k} \int^{\infty}_{t}   e^{-\mu z} dz \right) \cdot  e^{-\mu t} \\
&=& \frac{\lambda^k}{(\lambda + \mu)^k}  e^{-\mu t} \cdot  e^{-\mu t}   \\
&=& \left( \frac{\lambda}{\lambda + \mu} \right)^k e^{-2\mu t}  .
\end{eqnarray}
This completes the proof.
\end{proof}
\end{theorem}

\begin{figure}[ht!]
\hspace{-.5in}~\includegraphics[scale=.2]{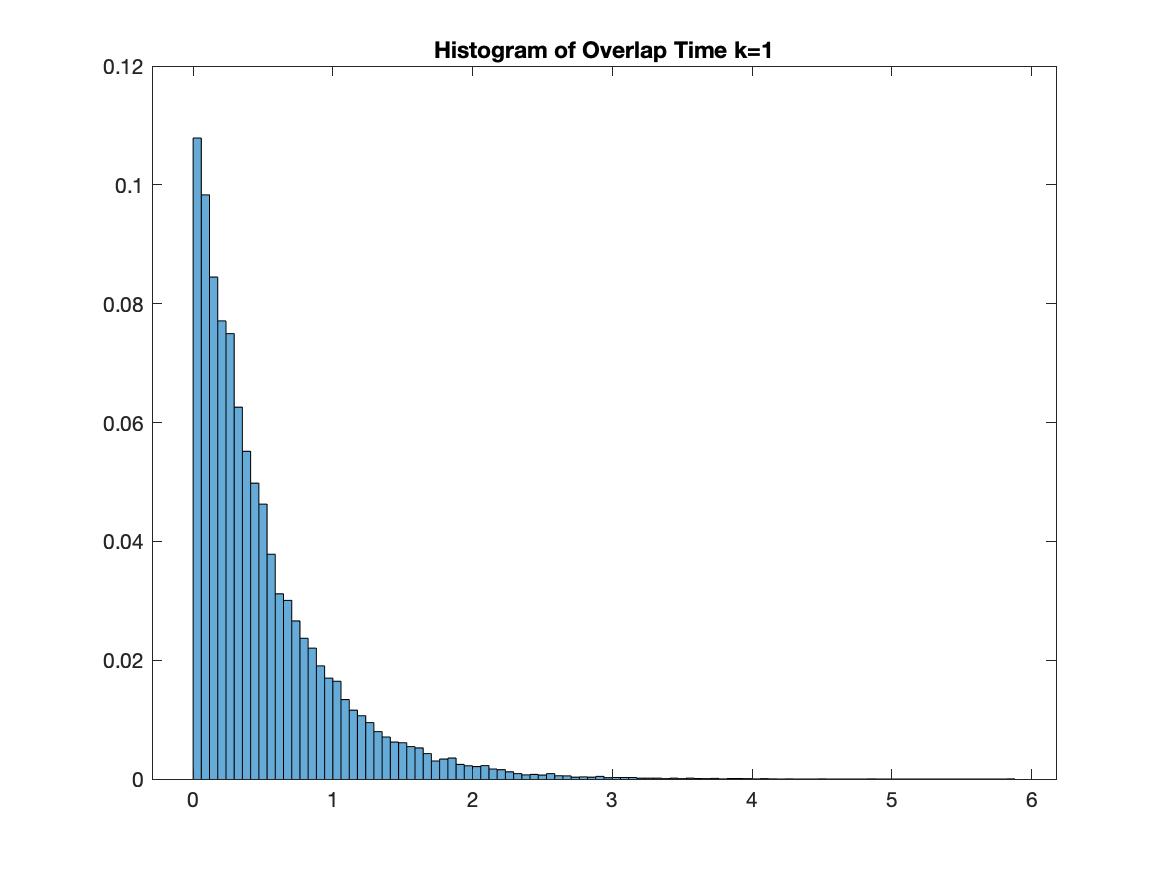}~\hspace{-.3in}~\includegraphics[scale =.2]{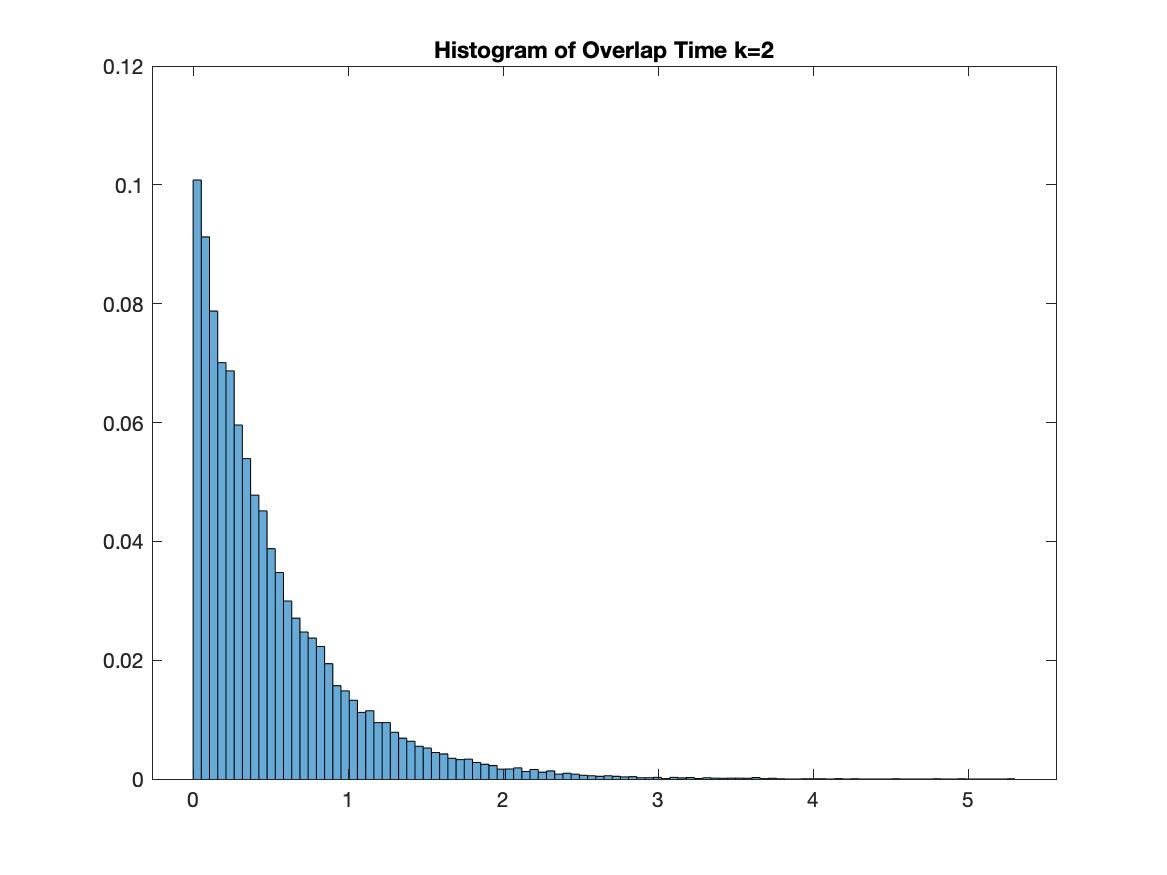}

\hspace{-.5in}~\includegraphics[scale=.2]{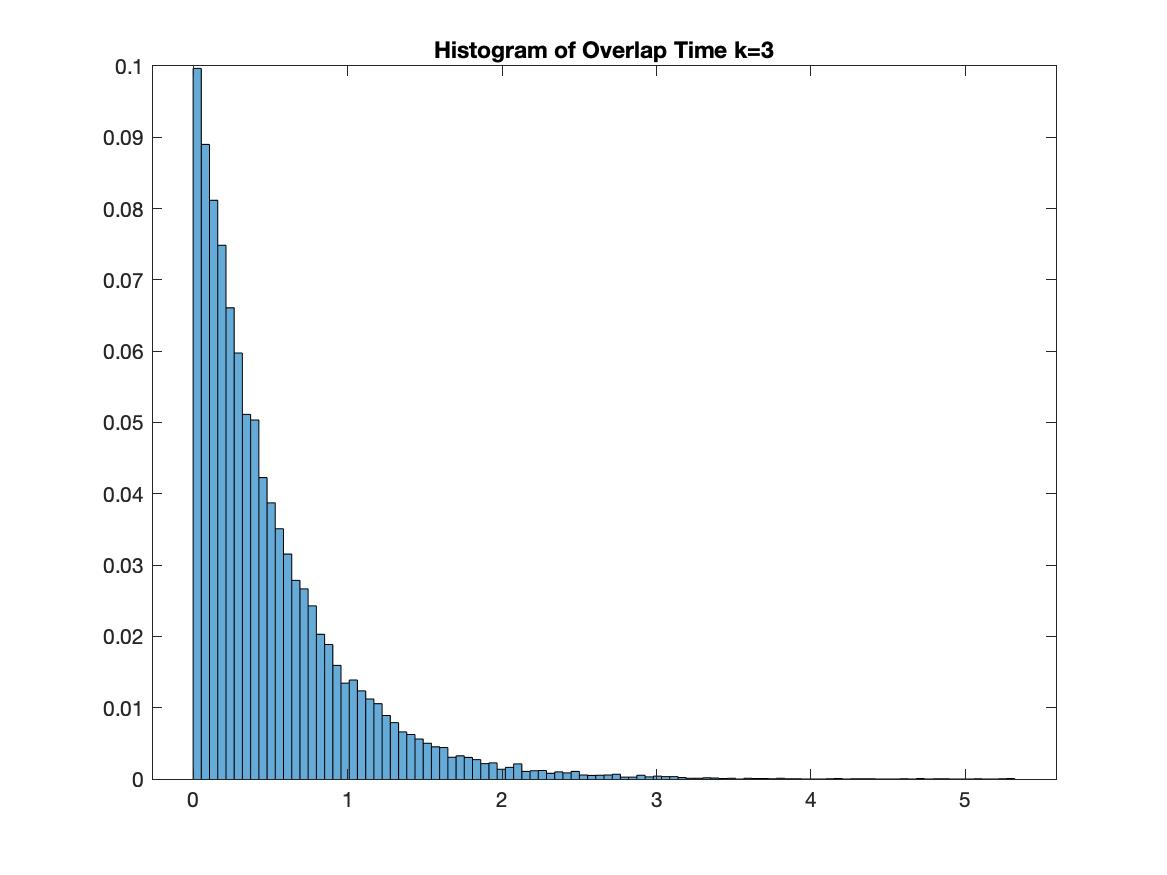}~\hspace{-.3in}~\includegraphics[scale =.2]{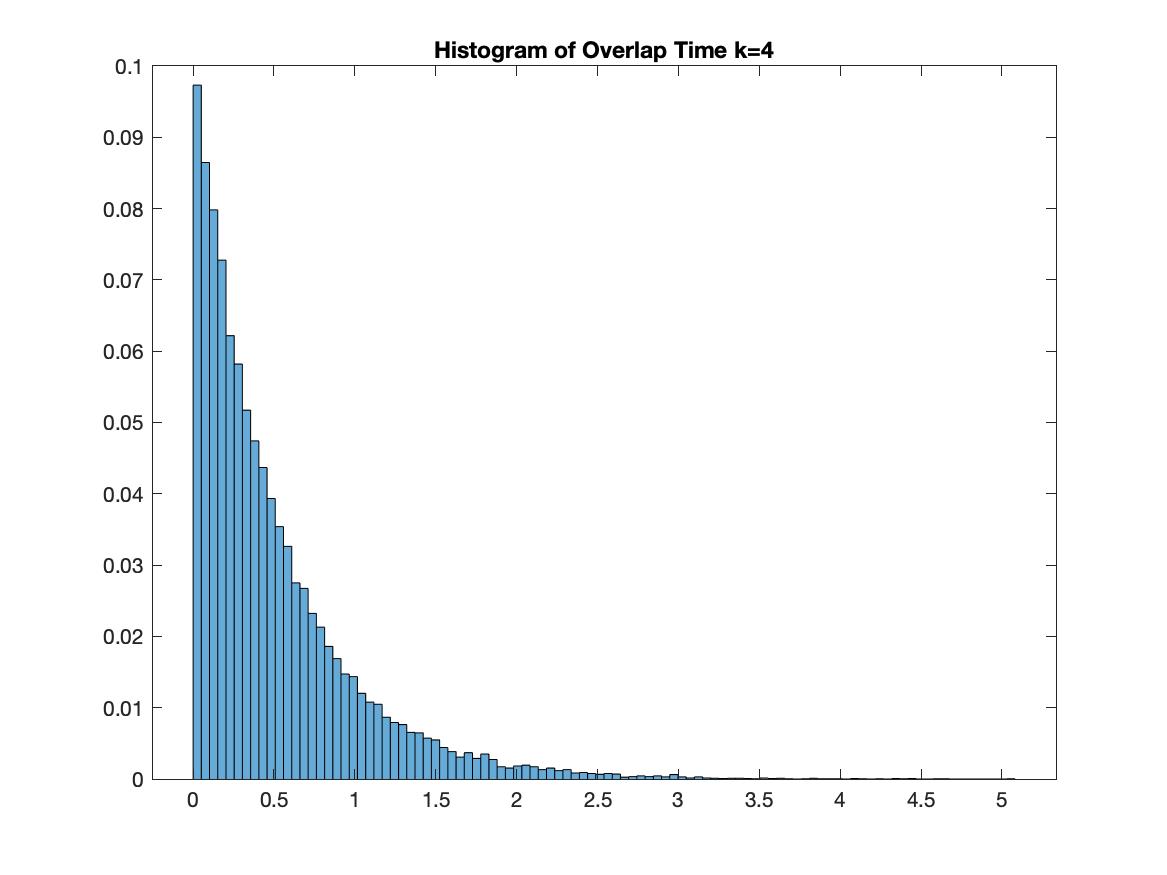}
 \captionsetup{justification=centering}
 \caption{Histograms for the Overlap Times Exponential Distributions $(\lambda = .5, \mu = 1)$.   Upper Left (k=1), Upper Right (k=2), Lower Left (k=3), Lower Right (k=4). }
 \label{fig:overlap_1}
\end{figure}

\begin{figure}[ht!]
\hspace{-.5in}~\includegraphics[scale=.2]{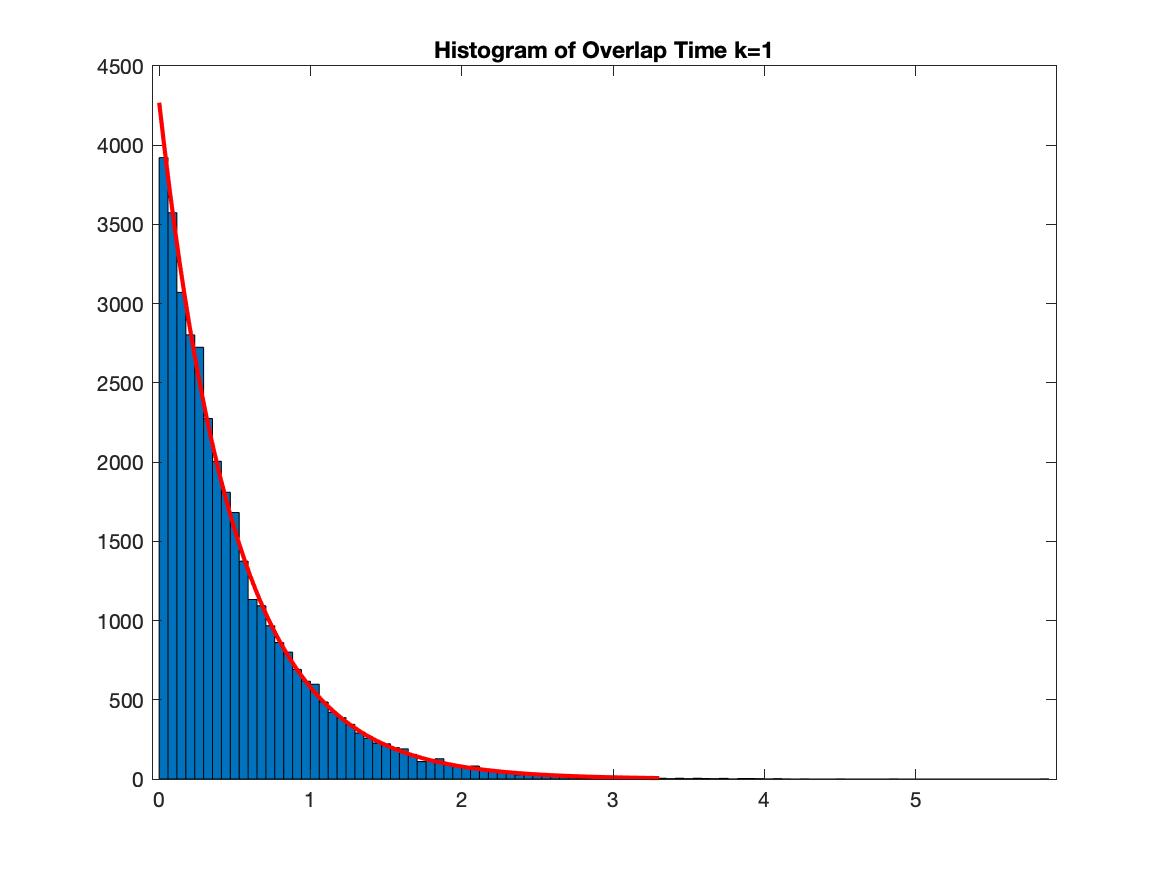}~\hspace{-.3in}~\includegraphics[scale =.2]{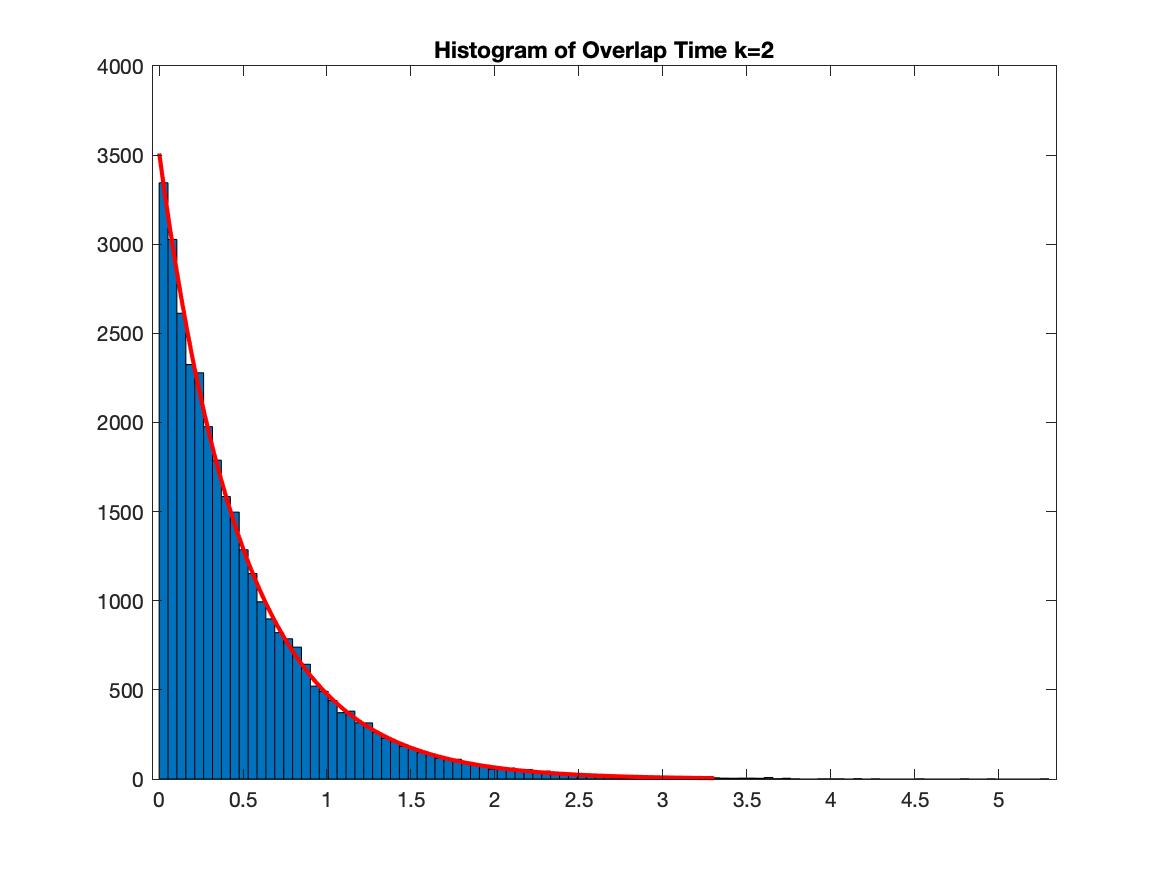}

\hspace{-.5in}~\includegraphics[scale=.2]{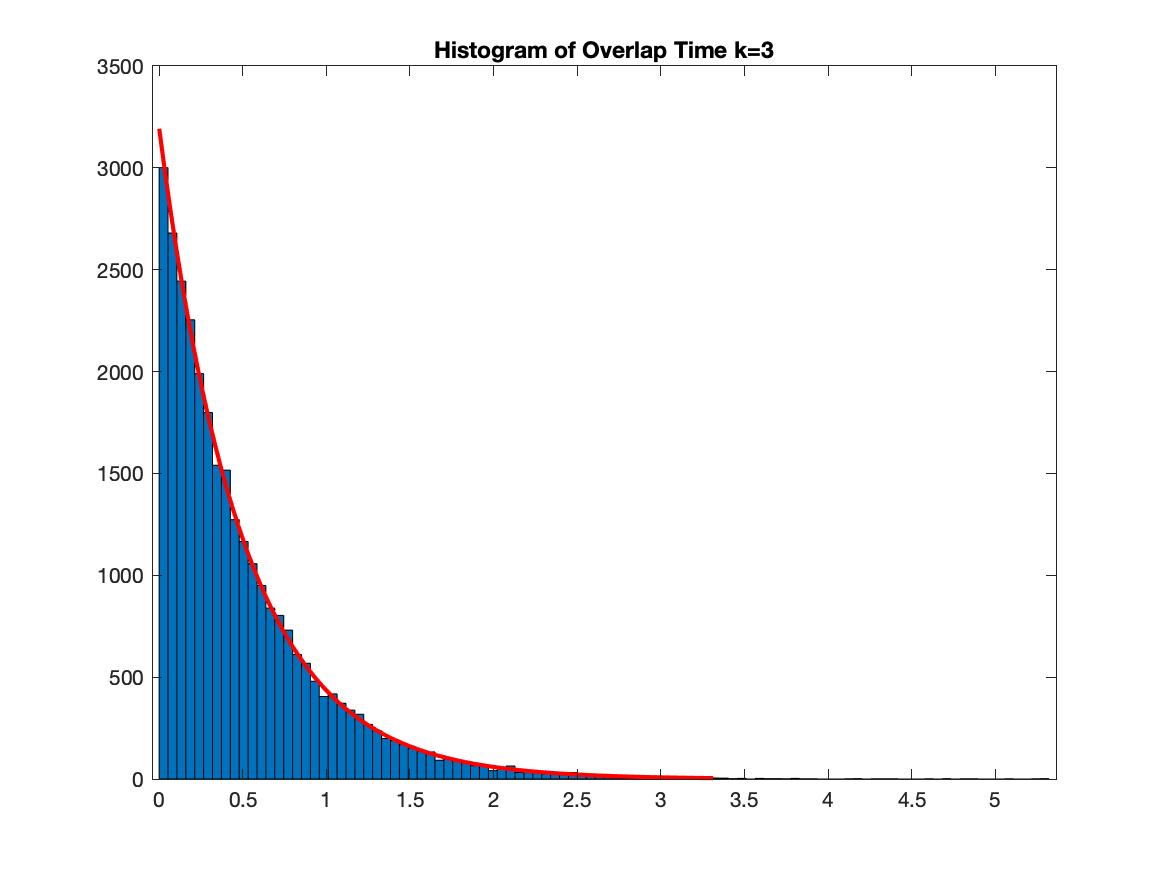}~\hspace{-.3in}~\includegraphics[scale =.2]{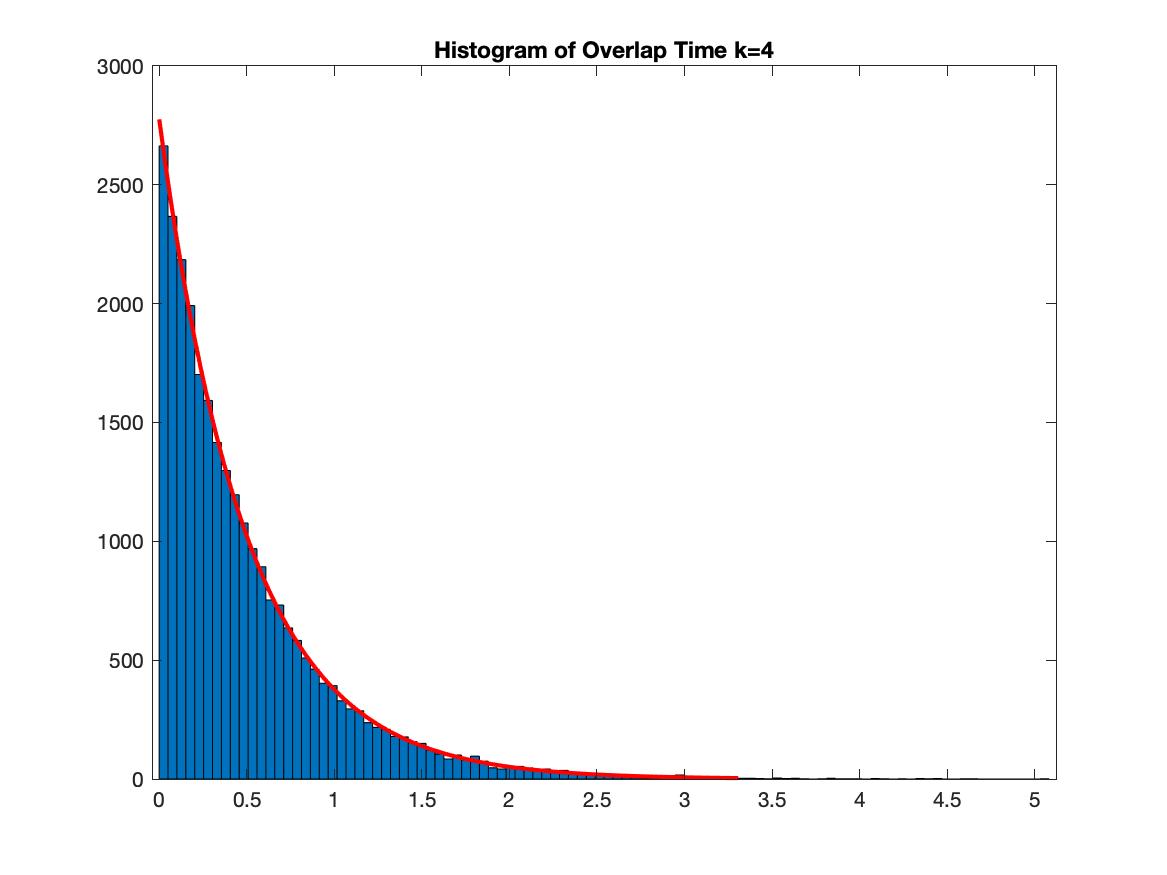}
 \captionsetup{justification=centering}
 \caption{Histograms for the Overlap Times Exponential Distributions $(\lambda = .5, \mu = 1)$.   Upper Left (k=1), Upper Right (k=2), Lower Left (k=3), Lower Right (k=4). }
 \label{fig:overlap_2}
\end{figure}

In Figure \ref{fig:overlap_1}, we plot the histograms of the overlap times of customers that have distance $k = \{ 1,2,3,4\}$ in the $M/M/\infty$ queue.  We find that the conditional distribution that the overlap time is greater than zero is an exponential distribution with rate $2\mu$.  However, why does the tail decay with rate $2\mu$?  One observation is that in order for the overlap time to be large, it must be that both service times are also large.  Since they each decay at rate $\mu$, then both being large is given by $2\mu$.  Note that this is different from the single server setting in \citet{palomooverlap} where the decay of overlap times is equal to $\mu$.  This is because the $(n+k)^{th}$ customer cannot leave until the $n^{th}$ customer has departed.  Thus, the decay only depends on the service time of the $n^{th}$ customer being large.

To see how well our theoretical analysis matches with our simulation experiments,  using maximum likelihood estimation, we fit an exponential distribution to the overlap times in Figure \ref{fig:overlap_2} showing that the exponential distribution is accurate as predicted by our theoretical analysis.  What is also apparent from Figure \ref{fig:overlap_2} is that the number of overlap times decreases as we increase the distance $k$.  This is consistent with the geometric-like decay of the coefficient in front of the exponential distribution given in Equation \ref{overlap_dist}.  Although we mostly have analyzed the $GI/GI/\infty$ queue in the context of Poisson arrivals, we can also find the steady state overlap time distribution for renewal arrivals when the service distribution is deterministic.  Theorem \ref{general_overlap}  below provides a proof of this result.




\begin{theorem} \label{general_overlap}
The steady state distribution of $O_{k}$ in the $GI/D/\infty$ queue is given by
\begin{eqnarray} \label{overlap_dist_det}
\mathbb{P} \left( O_{k} > t  \right) &=&  F^{(k)} \left( (\Delta - t)^+ \right)  .
\end{eqnarray}
\begin{proof}
\begin{eqnarray}
\mathbb{P} \left( O_{k} > t  \right) &=& \mathbb{P} \left( ( \Delta - \mathcal{A}_k  )^+ > t \right) \\
&=& \mathbb{P} \left( \Delta - \mathcal{A}_k  > t \right) \\
&=& \mathbb{P} \left( \mathcal{A}_k <  \Delta - t  \right) \\
&=& F^{(k)} \left( (\Delta - t)^+ \right)
\end{eqnarray}
where $ F^{(k)}(x)$ is the $k$-fold convolution of distribution function $F(x)$.  
\end{proof}
\end{theorem}

\begin{figure}[ht!]
\hspace{-.5in}~\includegraphics[scale=.2]{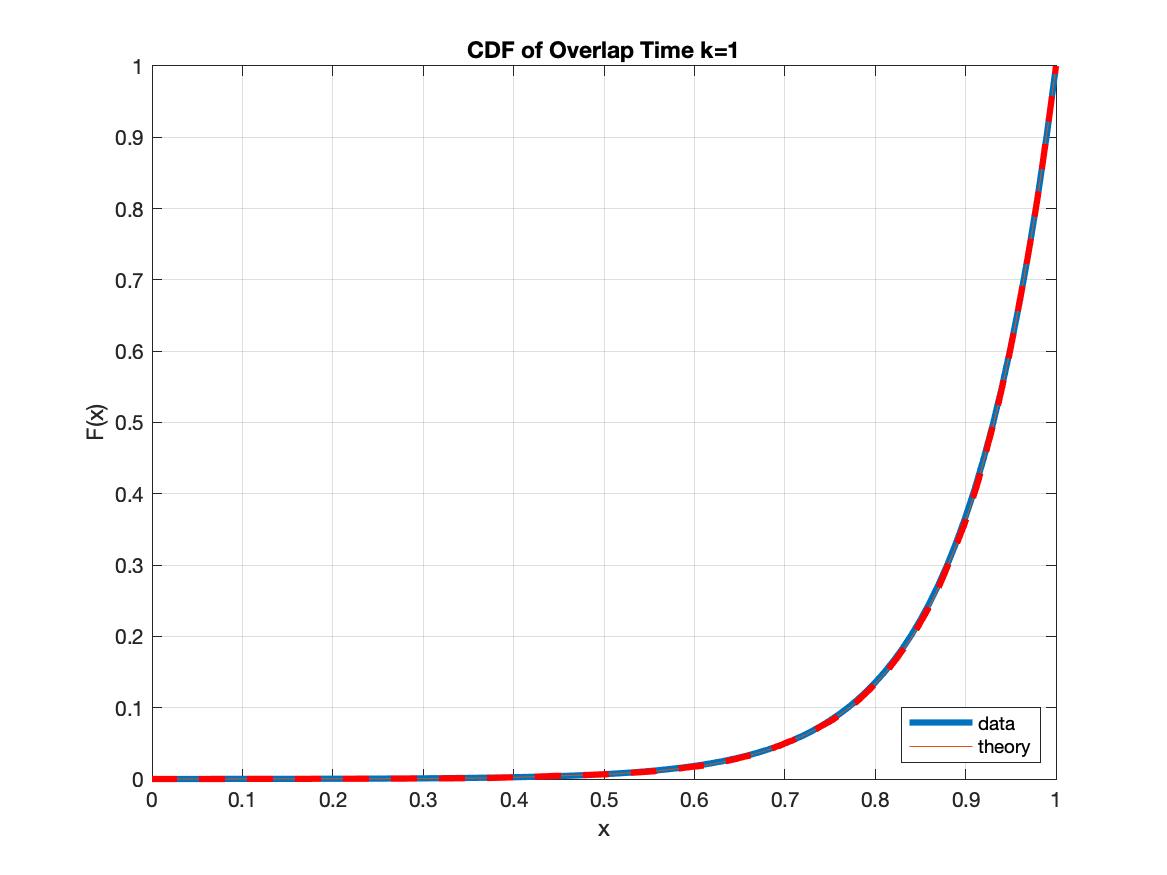}~\hspace{-.3in}~\includegraphics[scale =.2]{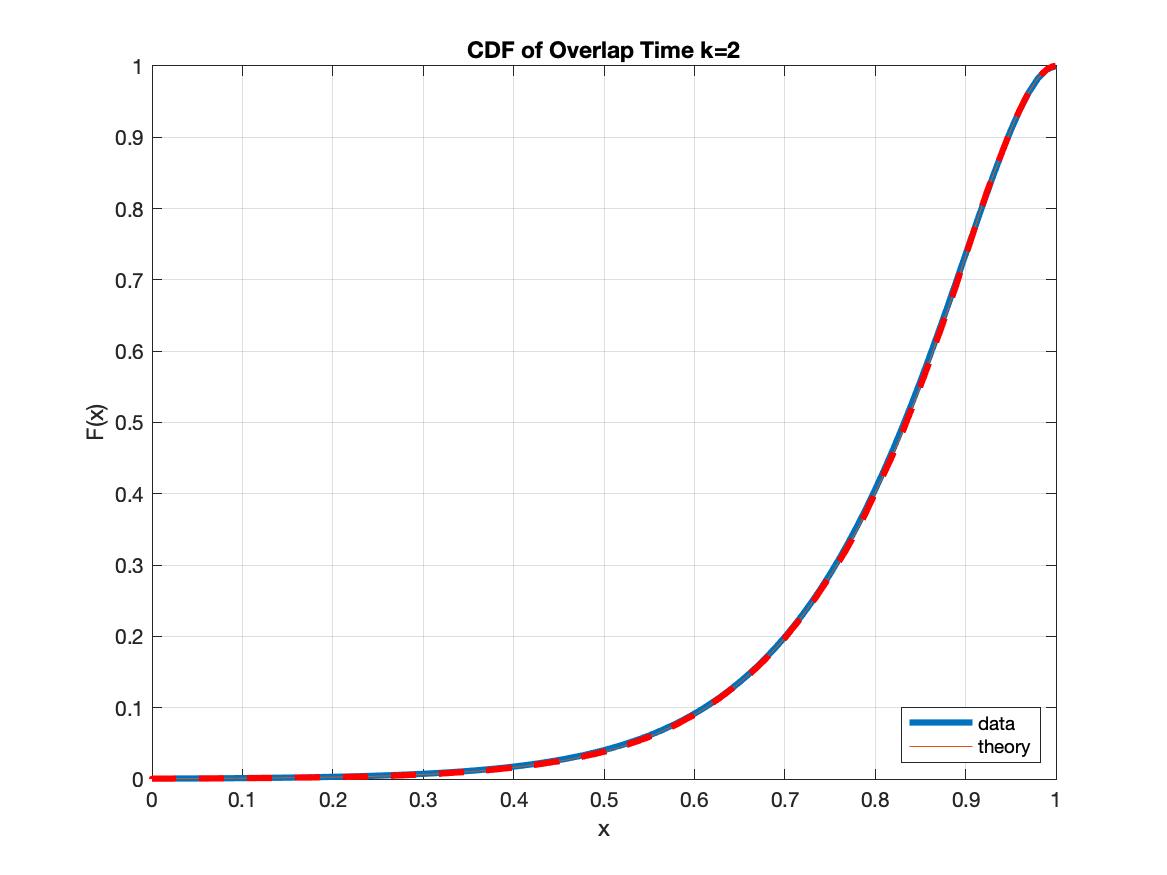}

\hspace{-.5in}~\includegraphics[scale=.2]{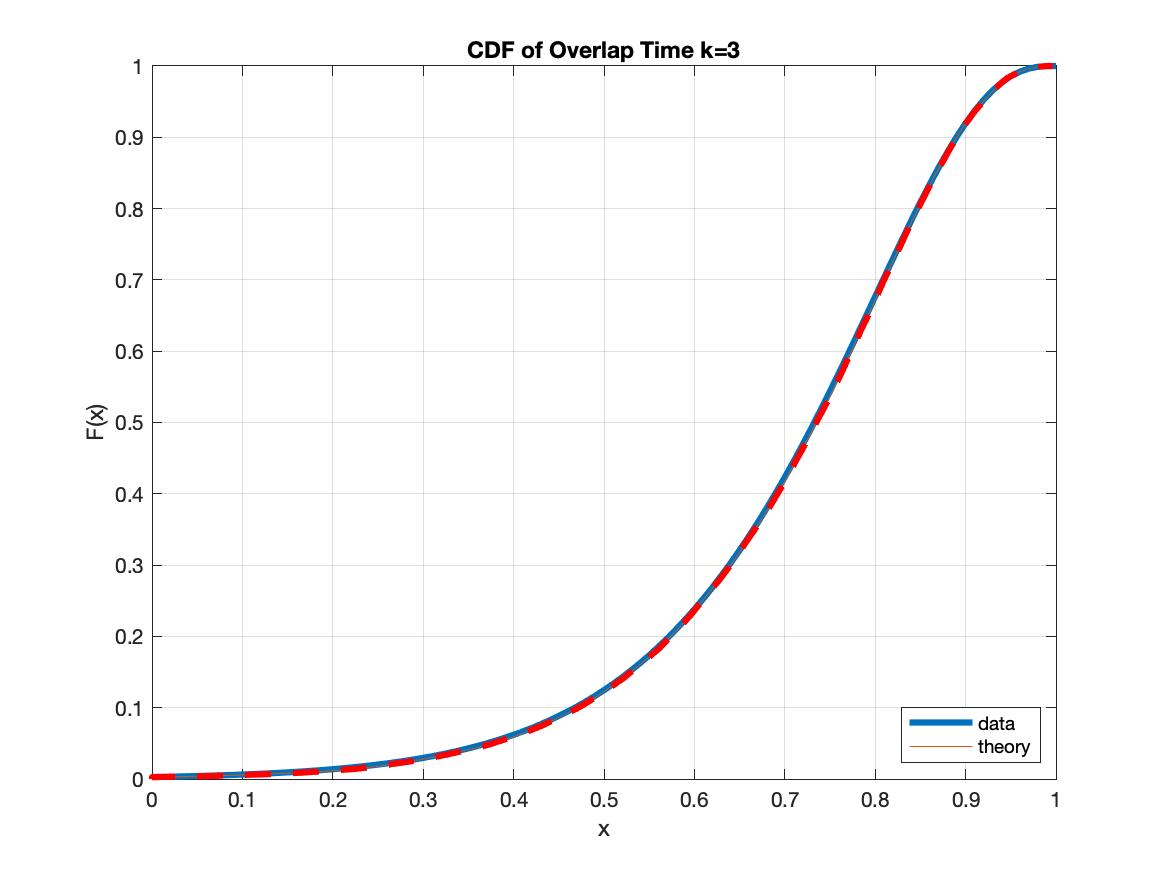}~\hspace{-.3in}~\includegraphics[scale =.2]{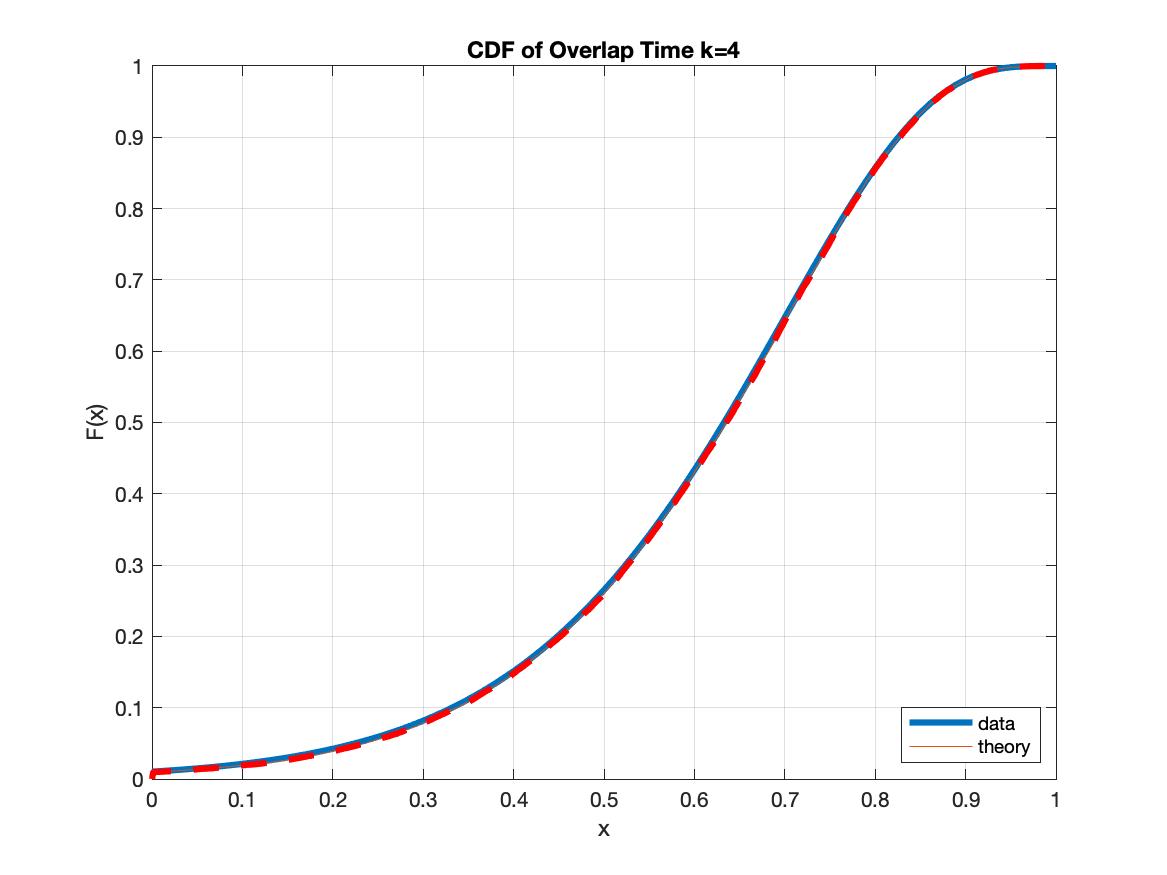}
 \captionsetup{justification=centering}
 \caption{CDF for the Overlap Times Deterministic Service $(\lambda = .5, \mu = 1)$.   Upper Left (k=1), Upper Right (k=2), Lower Left (k=3), Lower Right (k=4). }
 \label{fig:overlap_3}
\end{figure}



\begin{corollary}
In the $Gamma(\alpha,\lambda)/D/\infty$ queue, the steady state overlap distribution is given by
\begin{eqnarray}
\mathbb{P} \left( O_{k} > t  \right) &=&  \frac{\gamma(k \alpha, \lambda (\Delta - t)^+ )}{\Gamma(k \alpha)}.
\end{eqnarray}
\end{corollary}

In Figure \ref{fig:overlap_3}, we plot the cumulative distribution function (cdf) of the overlap times of the $M/D/\infty$  queue.  We observe that unlike the exponential service setting, the distribution changes significantly as we change the adjacency $k$. We find in all of the plots that the empirical distribution matches almost perfectly to the theoretical analysis given in Equation \ref{overlap_dist_det}.  It would be interesting to explore more distributions for the arrival and service random variables and see if we can compute them in closed form.  Suppose that we let the arrival and service distribution be Gamma distributed, then the resulting distribution can be described as a gamma difference distribution, see for example \citet{klar2015note}.  However, beyond this, it seems quite difficult to get results for other distributions in closed form.  

\section{Computing the Number of Overlaps} \label{Sec:Overlaps}

In addition to knowing how much time adjacent customers will spend together in a service system, it is important to also know how many customers one expects to overlap with.  In the context of epidemics, the more customers that one actually overlaps with increases the chances that one might contract the disease.  In this section, we restrict our analysis to that of the $M_t/G/\infty$ queue and leverage results from the $M_t/G/\infty$ queue in \citet{eick1993physics} to compute the exact distribution of overlapping customers.  In what follows, we assume without loss of generality that all queues start with zero customers.  Under this assumption, we know that the queue $M_t/G/\infty$ queue at time $t$ is given by the following expression
\begin{eqnarray}
Q^{\infty}(t) &=& \sum_{j=1}^{N^{(b)}(t)} \{A_j < t < A_j + S_j\}.
\end{eqnarray}
Thus, the number of people in the system at the time of arrival of the $k^{th}$ person is given by
\begin{eqnarray}
Q^{\infty}(A_{k}-) &=& \sum_{j=1}^{N^{(b)}(A_{k}-)} \{A_j < A_{k} < A_j + S_j\}. 
\end{eqnarray}
In order to compute the total number of overlapping customers, we also need to calculate the number of people that arrive during the service time of the $k^{th}$ customer as well.  The number of customers that arrive during the $k^{th}$ customer's service time is equal to 
\begin{eqnarray}
N^{(b)}(A_{k} + S_k) - N^{(b)}(A_{k}) . 
\end{eqnarray}
Adding the number of people upon arrival and the number that arrive during service, we arrive at an expression for $T_k$, the total number of people that the $k^{th}$ person will overlap with.  Thus, $T_k$ is equal to 
\begin{eqnarray}
T_k &=& N^{(b)}(A_{k} + S_k) - N^{(b)}(A_{k})  + \sum_{j=1}^{N^{(b)}(A_{k})} \{A_j < A_{k} < A_j + s_j\} .
\end{eqnarray}
It is important to note that the total number of overlapping customers for the $k^{th}$ arrival can be only computed when the $k^{th}$ arrival departs the queue and it is not known at the time of arrival.  Thus, this representation for the number of overlaps provides a methodology for computing the number of overlaps for each customer via simulation.  However, this representation is customer centered and not time centered.  In what follows, we provide a time centered perspective of the number of overlaps.

If a person arrives at time $t$, then we define $O(t)$ as the number of people that the arriving customer will overlap with.  $O(t)$ has the following expression
\begin{eqnarray}
O(t) &=& \underbrace{N^{(b)}(t + \mathcal{S}) - N^{(b)}(t)}_{\text{\# of customer arrivals during service}}  + \underbrace{\sum_{j=1}^{N^{(b)}(t)} \{A_j < t < A_j + S_j\}}_{\text{\# of customers upon arrival}} \\
&=& \underbrace{N^{(b)}(t + \mathcal{S}) - N^{(b)}(t)}_{\text{\# of customer arrivals during service}}   + \underbrace{Q^{\infty}(t)}_{\text{queue length at time t}} .
\end{eqnarray}
It is important to note here that the number of overlapping customers $O(t)$ is not completely known at time $t$.  It is only fully known after the service of the customer that arrives at time $t$. However, we can use the expression of $O(t)$ to compute the distribution of the number of customers that an arrival will overlap with at time $t$.  In this way, it is time centered as the expression depends on time and not the specific number of the customer arriving.  In addition to knowing the exact number of customers an arrival at time $t$ would overlap with or an exact sample path, we can also describe the exact distribution of the number of customers that an arrival at time $t$ would overlap.  One important ingredient to describing the exact distribution is knowing that the number of customers upon arrival is independent from the number of arrivals that arrive during service.  We provide the exact distribution in the following theorem.  

\begin{theorem}
The distribution of the number of overlapping customers in the $M_t/G/\infty$ queue at time $t$ is equal to 
\begin{eqnarray}
O(t) &\stackrel{D}{=}& \text{Poisson} \left( \int^{t+\mathcal{S}}_{t}\lambda(s)ds  +  \int^{t}_{0} \lambda(u) \overline{G}(t-u) du \right) .
\end{eqnarray}
Moreover, the transient state probabilities are given by the following integral expression
\begin{eqnarray}
\mathbb{P} \left( O(t) = k \right) &=& \int^{\infty}_{0} e^{ -(\Lambda(t,s) + q^{\infty}(t))} \frac{(\Lambda(t,s) + q^{\infty}(t))^k}{k!} dG(s) 
\end{eqnarray}
where $\Lambda(t,s)$ is defined as $\Lambda(t,s) = \int^{t+s}_t \lambda(u)du$.
\begin{proof}
\begin{eqnarray}
O(t) &=&  N^{(b)}(t + \mathcal{S}) - N^{(b)}(t)  + Q^{\infty}(t) \\
&\stackrel{D}
     {=}& \text{Poisson} \left( \int^{t+\mathcal{S}}_{t}\lambda(u) du  \right) + \text{Poisson} \left( q^{\infty}(t) \right) \\
&\stackrel{D}
     {=}& \text{Poisson} \left( \int^{t+\mathcal{S}}_{t}\lambda(u) du   \right) + \text{Poisson} \left(  \int^{t}_{0} \lambda(u) \overline{G}(t-u) du \right) \\
&\stackrel{D}
     {=}& \text{Poisson} \left( \int^{t+\mathcal{S}}_{t}\lambda(u) du  +  \int^{t}_{0} \lambda(u) \overline{G}(t-u) du  \right) .
\end{eqnarray}
This completes the proof.  
\end{proof}
\end{theorem}

\subsection{Number of Arrivals During Service}

It is important to remark that the distribution of the number of overlapping customers is not necessarily Poisson unless the service distribution is given by a discrete random variable.  To see this, note that even in the constant arrival rate case, an exponential service time yields a geometric random variable.  Moreover, a gamma distributed service time yields a negative binomial random variable.  In what follows, we calculate the number of arrivals during service for several service time distributions that are common in applied probability.  We also find an expression in terms of an expectation for any non-negative continuous random variable as well.  

\begin{proposition}
Let $\mathcal{S}$ be a Erlang$(\alpha, \mu)$ random variable, then the number of customers that arrive during service has a negative binomial distribution i.e.
\begin{eqnarray}
\mathbb{P} \left( N(t+\mathcal{S}) - N(t) = k \right) = {{k+\alpha-1}\choose{k} } ( 1 - \rho)^{\alpha} \rho^ k 
\end{eqnarray}
where $\rho = \frac{\lambda}{\mu + \lambda}$.
\begin{proof}
Since $\mathcal{S}$ is Erlang$(\alpha, \mu)$ distributed, we can decompose it into $\alpha$ i.i.d exponential random variables.  Thus, it suffices to prove the result and show the distribution is Geometric$\left( \rho \right)$.  
Thus, we need to show that if $\mathcal{S}$ is Exp($\mu$), then
\begin{eqnarray}
\mathbb{P} \left( N(t+\mathcal{S}) - N(t) = k \right) = ( 1 - \rho) \rho^ k .
\end{eqnarray}

Thus, we have
\begin{eqnarray}
\mathbb{P} \left( N(t+\mathcal{S}) - N(t) = k \right) &=& \int^{\infty}_0 e^{-\lambda x} \frac{(\lambda x)^k}{k!} \mu e^{-\mu x} dx \\
&=& \frac{\lambda^k \mu}{k!} \int^{\infty}_0 x^k e^{-(\lambda+\mu) x}  dx \\
&=& \frac{\lambda^k \mu}{(\lambda + \mu) k!} \int^{\infty}_0 x^k \cdot (\lambda + \mu) e^{-(\lambda+\mu) x}  dx \\
&=& \frac{\lambda^k \mu}{(\lambda + \mu) k!} \frac{k!}{(\lambda + \mu)^k} \\
&=& \frac{\lambda^k \mu}{(\lambda + \mu)^{k+1} }  \\
&=& \left( \frac{\lambda }{\lambda + \mu } \right)^k \cdot \left( \frac{\mu}{\lambda + \mu } \right) \\
&=& \left( \frac{\lambda }{\lambda + \mu } \right)^k \cdot \left( 1 - \frac{\lambda}{\lambda + \mu } \right) \\
&=& ( 1 - \rho) \rho^ k .
\end{eqnarray}
This completes the proof. 
\end{proof}
\end{proposition}

\begin{figure}[ht!]
\hspace{-.5in}~\includegraphics[scale=.2]{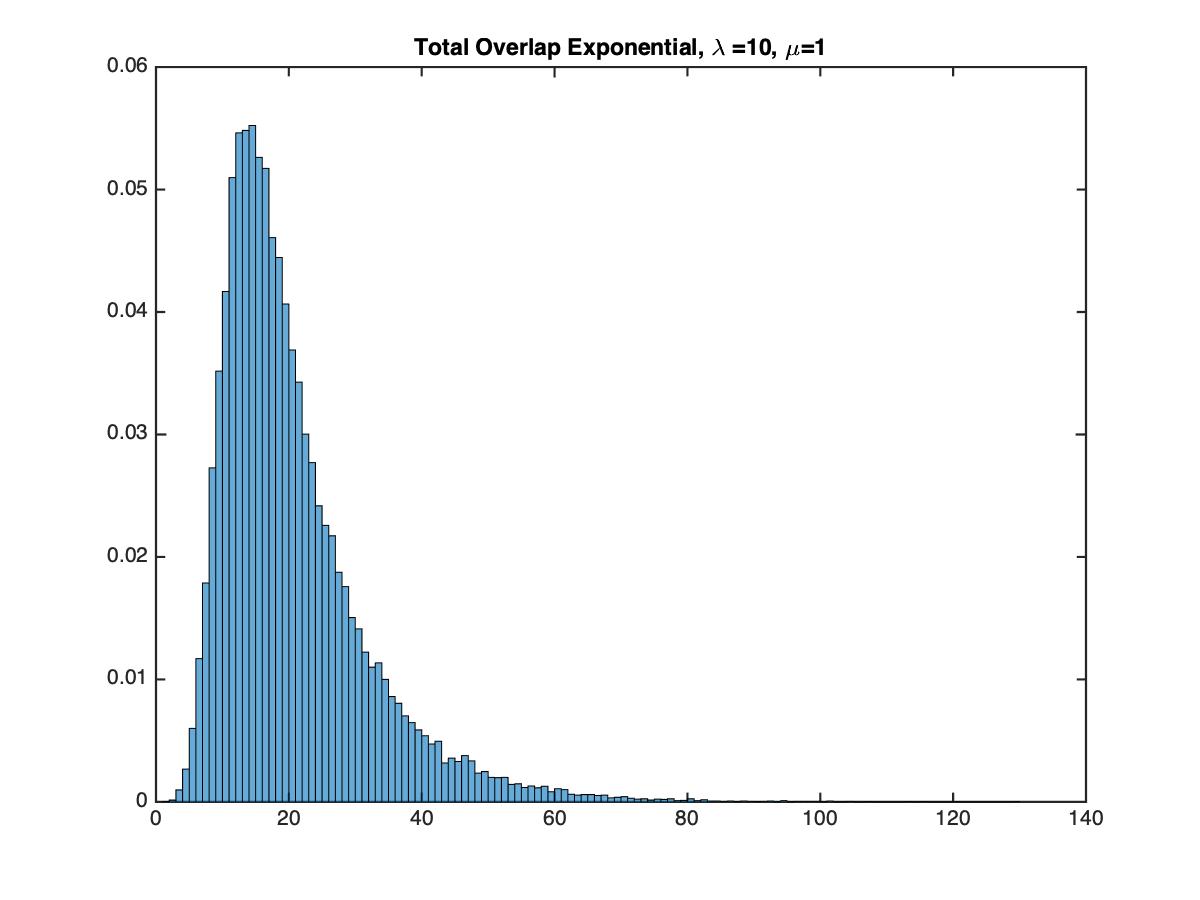}~\hspace{-.25in}~\includegraphics[scale=.2]{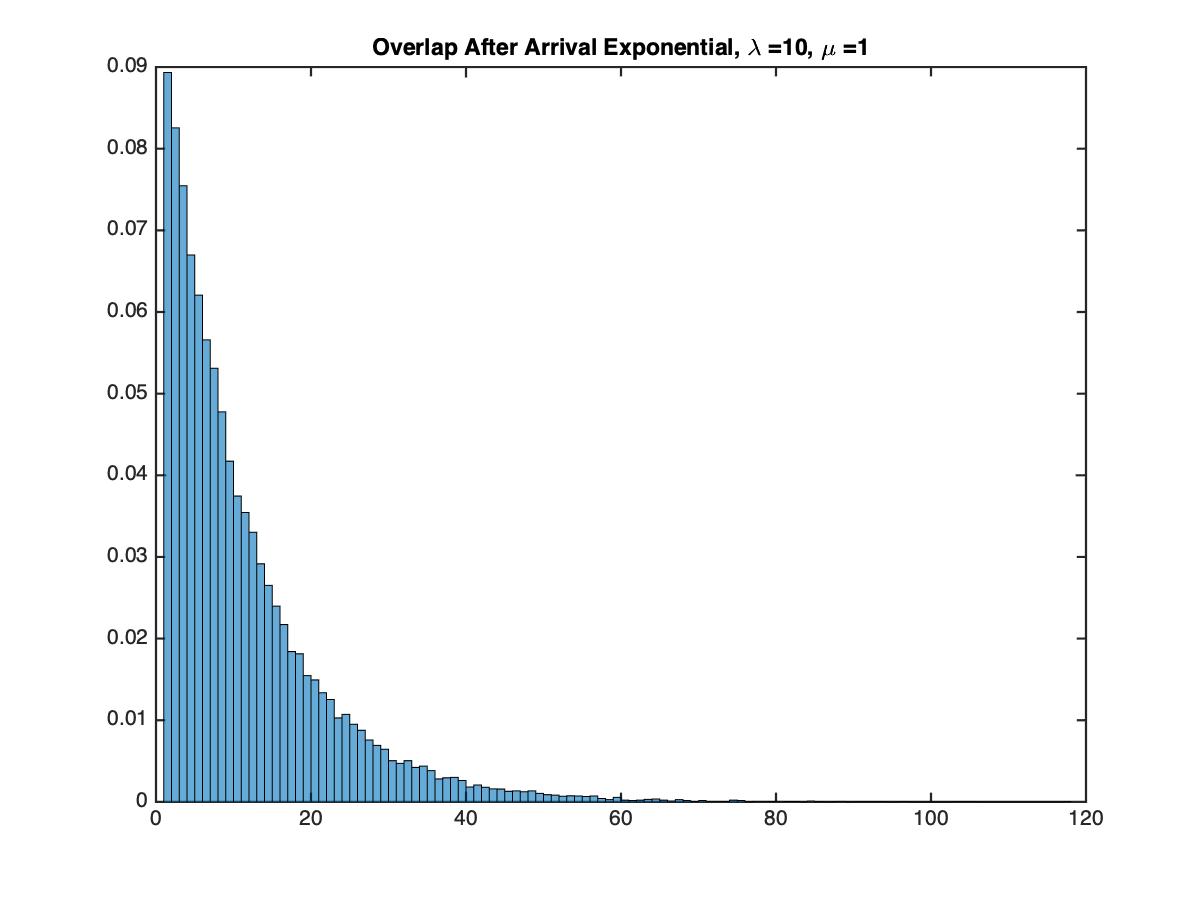}

\hspace{1in}~\includegraphics[scale=.2]{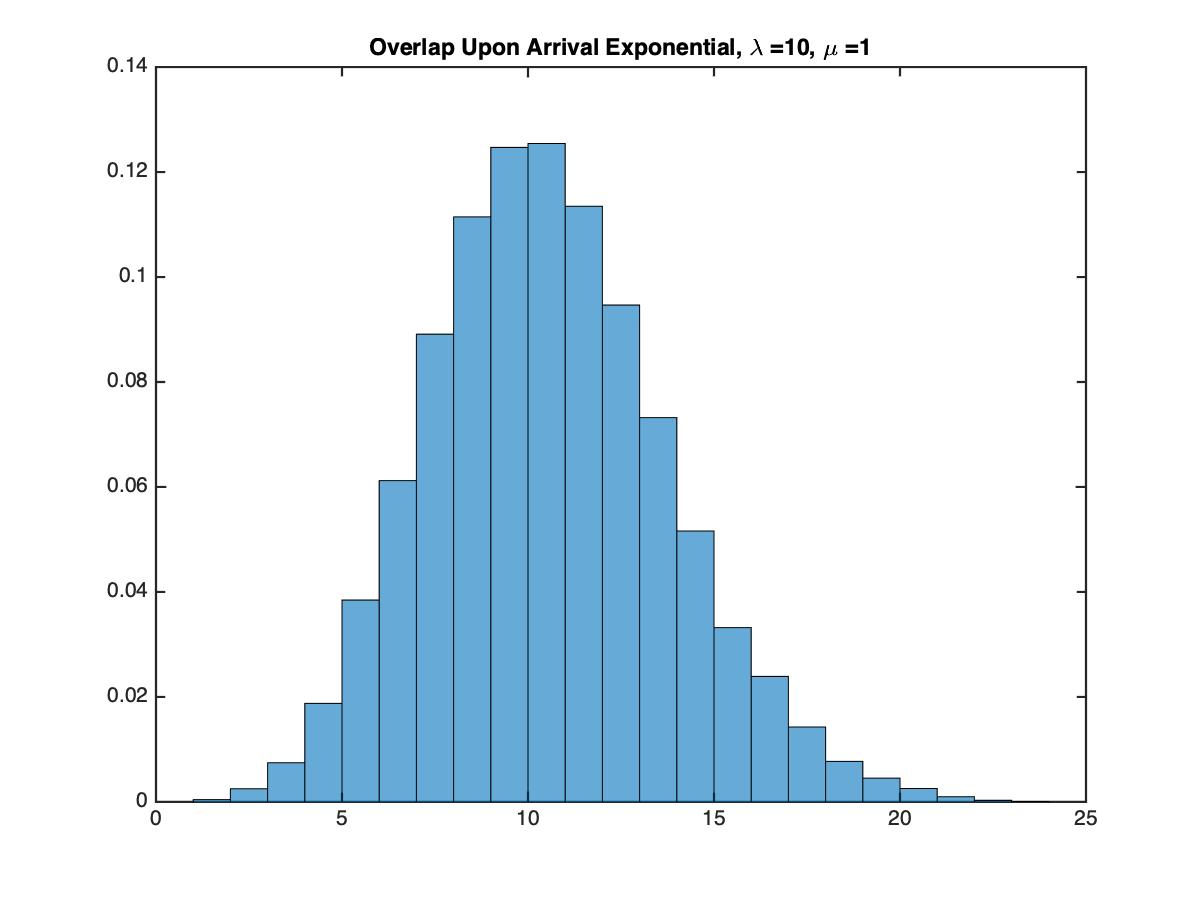}
 \captionsetup{justification=centering}
 \caption{Total Overlap (Exponential), $(\lambda = 10, \mu = 1)$ (Top Left). \\ 
 Overlap After Arrival (Exponential), $(\lambda = 10, \mu = 1)$ (Top Right). \\
 Overlap Upon Arrival (Exponential), $(\lambda = 10, \mu = 1)$ (Bottom). }
 \label{fig:overlap_all_8}
\end{figure}

In Figure \ref{fig:overlap_all_8}, we plot the distribution for the total number of overlaps, the number of overlaps upon arrival and the number overlaps during service.  We see that the distribution of the overlaps upon arrival is still Poisson and the number of overlaps during service is geometric as predicted by our theoretical results.  The total number is a the convolution of the two as well.  Thus, we have a good sense of how many customers might overlap with during their service experience.  

\begin{proposition}
Let $\mathcal{S}$ be a Uniform$(a, b)$ random variable, then the number of customers that arrive during service has the following distribution
\begin{eqnarray}
\mathbb{P} \left( N(t+\mathcal{S}) - N(t) = k \right) &=& \frac{\gamma(k+1,\lambda b) - \gamma(k+1,\lambda  a) }{\lambda (b-a) k!}.
\end{eqnarray}

\begin{proof}
\begin{eqnarray}
\mathbb{P} \left( N(t+\mathcal{S}) - N(t) = k \right) &=& \int^{b}_a e^{-\lambda x} \frac{(\lambda x)^k}{k!} \cdot \frac{1}{b-a}  dx \\
&=& \frac{\lambda^k}{(b-a) k!}\int^{b}_a x^k e^{-\lambda x}  dx \\
&=& \frac{\lambda^k}{(b-a) k!} \left( \int^{b}_0 x^k e^{-\lambda x}  dx - \int^{a}_0 x^k e^{-\lambda x}  dx \right) \\
&=& \frac{\gamma(k+1,\lambda b) - \gamma(k+1,\lambda  a) }{\lambda (b-a) k!}.
\end{eqnarray}
This completes the proof. 
\end{proof}
\end{proposition}

\begin{proposition}
Let $\mathcal{S}$ be a Truncated Normal$(a, b, \mu, \sigma^2)$ random variable, then the number of customers that arrive during service has the following distribution
\begin{eqnarray}
\mathbb{P} \left( N(t+\mathcal{S}) - N(t) = k \right) &=& \frac{\gamma(k+1,\lambda b) - \gamma(k+1,\lambda  a) }{\lambda (b-a) k!}.
\end{eqnarray}

\begin{proof}
\begin{eqnarray}
\mathbb{P} \left( N(t+\mathcal{S}) - N(t) = k \right) &=& \int^{b}_a e^{-\lambda x} \frac{(\lambda x)^k}{k!} \cdot \frac{e^{-(x-\mu)^2/(2\sigma^2)}}{\sqrt{2\pi} \sigma \left( \Phi(b)-\Phi(a) \right)}  dx \\
&=& \frac{\lambda^k}{\sqrt{2\pi} \sigma (\Phi(b)-\Phi(a)) k!}\int^{b}_a x^k e^{-\lambda x}  e^{-(x-\mu)^2/(2\sigma^2)}  dx \\
&=& \frac{\lambda^k e^{\lambda \sigma^2/2 - \mu }}{ \sqrt{2\pi} \sigma (\Phi(b)-\Phi(a)) k!}\int^{b}_a x^k  e^{-(x- \left( \mu - \lambda \sigma^2 \right))^2/(2\sigma^2)}  dx \\
&=&  \frac{\lambda^k e^{\lambda \sigma^2/2 - \mu }}{ k!}\int^{b}_a x^k  \frac{ e^{-(x- \left( \mu - \lambda \sigma^2 \right))^2/(2\sigma^2)} }{\sqrt{2\pi} \sigma (\Phi(b)-\Phi(a))}  dx \\
&=&  \frac{\lambda^k e^{\lambda \sigma^2/2 - \mu }}{ k!} M(k,a,b,\mu - \lambda \sigma^2/2, \sigma^2)
\end{eqnarray}
where $M(k,a,b,\mu - \lambda \sigma^2/2, \sigma^2)$ is the $k^{th}$ moment of the truncated normal distribution and its expression is given in Theorem 2.3 of \citet{pender2015truncated}.  Note that this can also be viewed as a change of measure of a Gaussian distribution and taking the truncated moment.  
\end{proof}
\end{proposition}

\begin{proposition} \label{mix_geo}
Let $\mathcal{S}$ be a Hyper-exponential$(\vec{p}, \vec{\mu}, \ell)$ random variable, then the number of customers that arrive during service has the following distribution 
\begin{eqnarray}
\mathbb{P} \left( N(t+\mathcal{S}) - N(t) = k \right) = \sum^{\ell}_{j=1} p_j ( 1 - \rho_j) \rho_j^ k 
\end{eqnarray}
where $\rho_j = \frac{\lambda}{\mu_j + \lambda}$.
\begin{proof}

Thus, we have
\begin{eqnarray}
\mathbb{P} \left( N(t+\mathcal{S}) - N(t) = k \right) &=& \sum^{\ell}_{j=1} \int^{\infty}_0 e^{-\lambda x} \frac{(\lambda x)^k}{k!} p_j\mu_j e^{-\mu_j x} dx \\
&=& \sum^{\ell}_{j=1} p_j ( 1 - \rho_j) \rho_j^k.
\end{eqnarray}
\end{proof}
\end{proposition}

From Proposition \ref{mix_geo}, we observe that a hyper-exponential service time turns the number of arrival after service into a convex combination of geometric probabilities.  This is nice from an analytical point of view since the distribution is a weighted sum of geometric distributions and the associated random variable is a Mixed-Geometric random variable, which is the discrete analogue of the hyper-exponential distribution. Note that the Mixed-Geometric is not called hyper-geometric as hyper-geometric refers to another distribution already.   

In addition to computing the distribution of arrivals during service for different service distribution, we can also compute the distribution in terms of a Laplace transform times a moment of the service distribution.  This is given in Theorem \ref{general_overlap} below.  

\begin{proposition} \label{general_overlap_2}
Let $\mathcal{S}$ be a continuous non-negative random variable, then the number of customers that arrive during service has the following distribution 
\begin{eqnarray}
\mathbb{P} \left( N(t+\mathcal{S}) - N(t) = k \right) = \frac{\lambda^k}{k!} \mathbb{E} \left[ \mathcal{S}^k e^{-\lambda \mathcal{S}} \right] = \frac{\lambda^k}{k!} \mathbb{E} \left[ e^{-\lambda\mathcal{S} + k \log(\mathcal{S}) } \right] .
\end{eqnarray}
\begin{proof}
We have the following expression
\begin{eqnarray}
\mathbb{P} \left( N(t+\mathcal{S}) - N(t) = k \right) &=&  \int^{\infty}_0 e^{-\lambda x} \frac{(\lambda x)^k}{k!} g(x) dx \\
&=& \frac{\lambda^k}{k!} \int^{\infty}_0 e^{-\lambda x} x^k g(x) dx \\
&=& \frac{\lambda^k}{k!} \mathbb{E} \left[ \mathcal{S}^k e^{-\lambda \mathcal{S}} \right] \\
&=& \frac{\lambda^k}{k!} \mathbb{E} \left[ e^{-\lambda\mathcal{S} + k \log(\mathcal{S}) } \right] \\
&=& \frac{(-\lambda)^k}{k!} \frac{\partial^k}{\partial \lambda^k }\mathbb{E} \left[ e^{-\lambda\mathcal{S}  } \right] .
\end{eqnarray}
\end{proof}
\end{proposition}

\begin{proposition} \label{coro_d}
For the $M_t/D/\infty$ queue the overlap distribution at time $t$ is equal to 
\begin{eqnarray}
O(t) &\stackrel{D}{=}& \text{Poisson} \left(  \int^{t+\Delta}_{(t-\Delta)^+} \lambda(u) du \right) .
\end{eqnarray}
\begin{proof}
\begin{eqnarray}
O(t) &\stackrel{D}{=}& \text{Poisson} \left( \int^{t+\mathcal{S}}_{t}\lambda(u) du  +  \lambda \int^{t}_{0} \overline{G}(t-u) du  \right) \\
&\stackrel{D}{=}& \text{Poisson} \left( \int^{t+\Delta}_{t}\lambda(u) du  +   \int^{t}_{(t-\Delta)^+} \lambda(u) du  \right) \\
&\stackrel{D}{=}& \text{Poisson} \left(  \int^{t+\Delta}_{(t-\Delta)^+} \lambda(u) du \right) .
\end{eqnarray}
This completes the proof.
\end{proof}
\end{proposition}



\begin{figure}[ht!]
\hspace{-.5in}~\includegraphics[scale=.2]{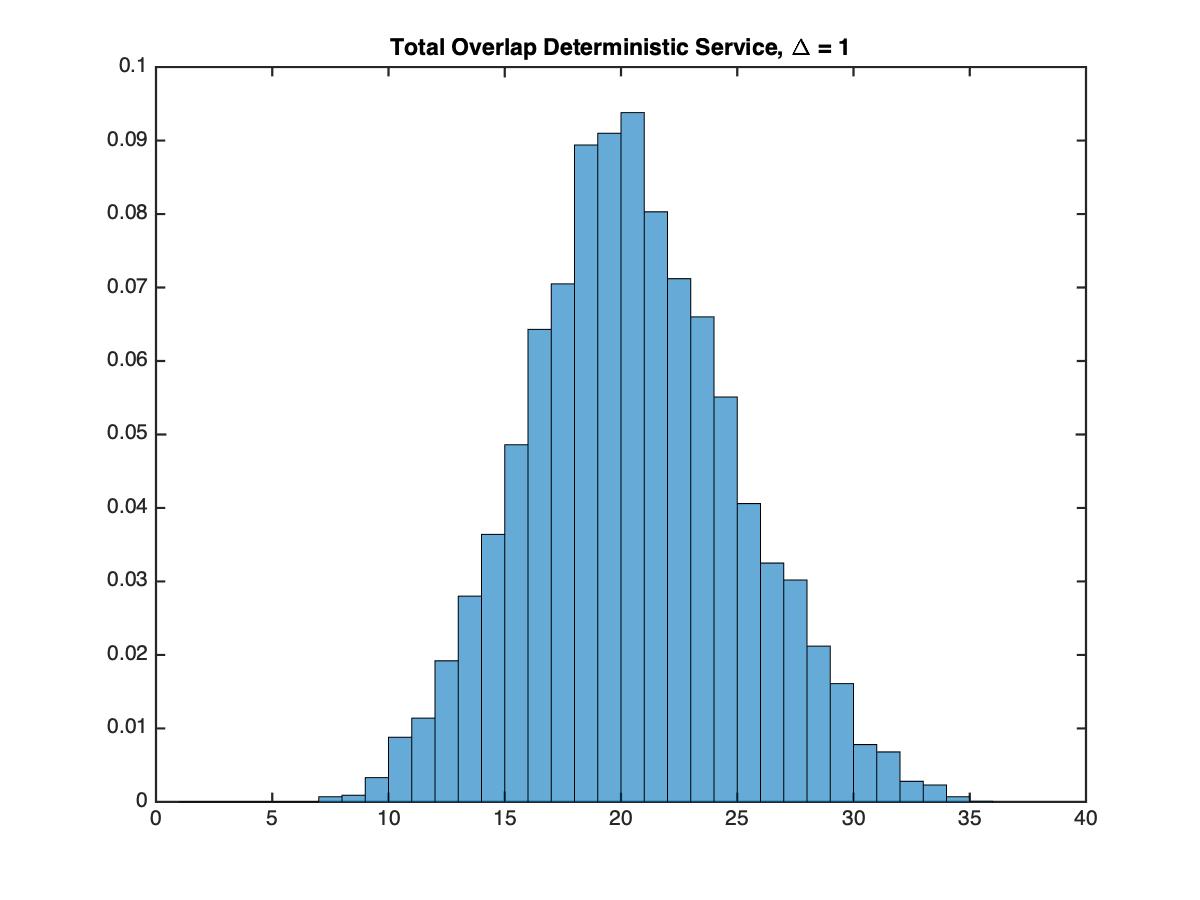}~\hspace{-.25in}~\includegraphics[scale=.2]{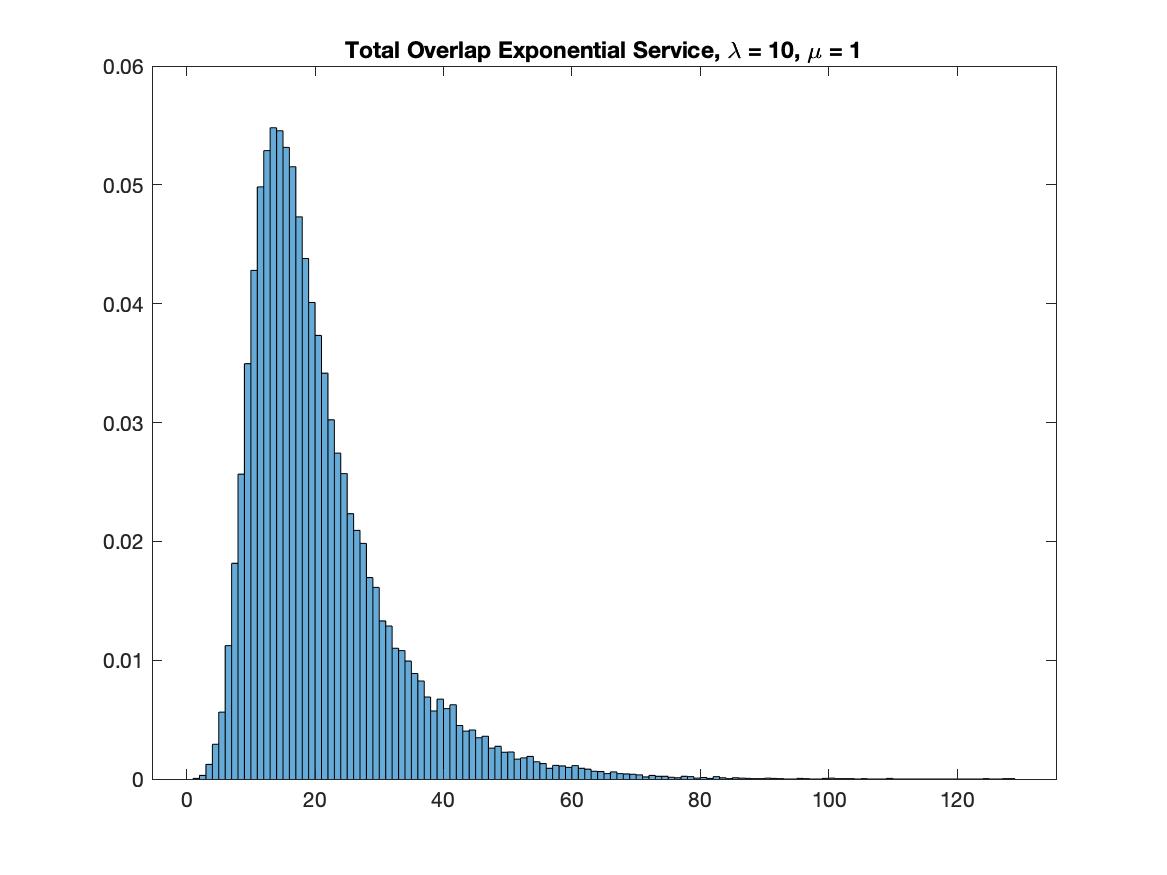}
 \captionsetup{justification=centering}
 \caption{Total Overlap Deterministic Service $(\lambda = 10, \Delta = 1)$ (Left). \\ 
 Total Overlap Exponential Service $(\lambda = 10, \mu = 1)$ (Right).}
 \label{fig:overlap_all_4}
\end{figure}

In Figure \ref{fig:overlap_all_4}, we plot the total overlap distribution for deterministic service and exponential service.  We see that while both overlap distributions have a mean of about 20, we see very different shapes in their distributions.  For the deterministic case, we know the overlap is Poisson, however, for the exponential setting it is the sum of a Poisson and geometric random variables.  The deterministic service setting results in a more symmetric distribution while the exponential setting is positively skewed.  

\begin{corollary}
For the $M_t/D(\ell)/\infty$ queue the overlap distribution at time $t$ is equal to 
\begin{eqnarray}
O(t) &\stackrel{D}{=}& \text{Poisson} \left( \sum^{\ell}_{j=1} \int^{t+\Delta_j}_{(t-\Delta_j)^+} p_j \lambda(u) du \right) .
\end{eqnarray}
\begin{proof}
The proof follows immediately from Poisson thinning and Proposition \ref{coro_d}
\end{proof}
\end{corollary}

\begin{theorem}
For the $M/M/\infty$ queue the overlap distribution at time $t$ is equal to 
\begin{eqnarray}
\mathbb{P} \left( O(t) = k \right) 
&=&  \rho^k (1-\rho)  e^{\frac{\mu}{\lambda} q^{\infty}(t)} \frac{\Gamma \left(k+1, \frac{\lambda + \mu}{\lambda} q^{\infty}(t) \right) }{\Gamma(k+1) }
\end{eqnarray}
and the mean and variance are given by 
\begin{eqnarray}
E[O(t)] 
&=&  \frac{\lambda}{\mu} ( 2 - e^{-\mu t} ),
\end{eqnarray}
\begin{eqnarray}
\mathrm{Var}[O(t)] &=&  \frac{\lambda^2}{\mu^2} + \frac{\lambda}{\mu} ( 1 - e^{-\mu t} ) .
\end{eqnarray}
Moreover, we have that $O(t)$ can be decomposed into a sum of geometric and Poisson random variables i.e.
\begin{eqnarray}
O(t) &\stackrel{D}{=}&  \mathrm{Geometric} \left( \frac{\lambda}{\lambda + \mu} \right) + \mathrm{Poisson}(q^{\infty}(t) ).
\end{eqnarray}
\begin{proof}
\begin{eqnarray}
\mathbb{P} \left( O(t) = k \right) &=& \int^{\infty}_{0} e^{ -(\lambda s + q^{\infty}(t))} \frac{(\lambda s + q^{\infty}(t))^k}{k!} \mu e^{-\mu s} ds \\
&=&\frac{\mu}{\lambda} e^{\frac{\mu}{\lambda} q^{\infty}(t)} \int^{\infty}_{ q^{\infty}(t)} e^{-u} \frac{ u^k}{k!} e^{-\frac{\mu}{\lambda} u } du, \quad  \quad u = \lambda s + q^{\infty}(t) \\
&=&\frac{\mu}{\lambda} e^{\frac{\mu}{\lambda} q^{\infty}(t)} \int^{\infty}_{ q^{\infty}(t)} \frac{ u^k}{k!} e^{-\frac{\lambda + \mu}{\lambda} u } du \\
&=&\frac{\mu}{\lambda} e^{\frac{\mu}{\lambda} q^{\infty}(t)} \frac{\Gamma \left(k+1, \frac{\lambda + \mu}{\lambda} q^{\infty}(t) \right) }{\Gamma(k+1) \cdot \left(\frac{\lambda + \mu}{\lambda} \right)^{k+1} } .
\end{eqnarray}
Finally for the mean, we have that 
\begin{eqnarray}
E[O(t)] 
&=&  E\left[ N^{(b)}(t + \mathcal{S}) - N^{(b)}(t) \right] + E[Q^{\infty}(t)] \\
&=& \lambda E[ \mathcal{S}] + \lambda \int^{t}_{0} \overline{G}(t-u) du \\
&=& \frac{\lambda}{\mu} + \frac{\lambda}{\mu} ( 1 - e^{-\mu t} )  \\
&=&   \frac{\lambda}{\mu} ( 2 - e^{-\mu t} ) 
\end{eqnarray}
Moreover, the variance can be calculated as follows
\begin{eqnarray}
\mathrm{Var}[O(t)] 
&=&  \mathrm{Var}\left[ N^{(b)}(t + \mathcal{S}) - N^{(b)}(t) \right] + \mathrm{Var}[Q^{\infty}(t)] \\
&=& \lambda E[\mathcal{S}] + \lambda^2 \mathrm{Var}[ \mathcal{S}] + \lambda \int^{t}_{0} \overline{G}(t-u) du \\
&=& \frac{\lambda^2}{\mu^2} + \frac{\lambda}{\mu} ( 2 - e^{-\mu t} )  .
\end{eqnarray}
This completes the proof.  
\end{proof}
\end{theorem}



\begin{proposition}
For the $M/H_\ell/\infty$ queue the overlap distribution at time $t$ is equal to 
\begin{eqnarray}
\mathbb{P} \left( O(t) = k \right) &=& \sum^\ell_{j=0} \rho^k_j (1-\rho_j)  e^{\frac{\mu_j}{\lambda} q^{\infty}_j(t)} \frac{\Gamma \left(k+1, \frac{\lambda + \mu_j}{\lambda} q^{\infty}_j(t) \right) }{\Gamma(k+1) }
\end{eqnarray}
and the mean is given by 
\begin{eqnarray}
E[O(t)] 
&=& \sum^\ell_{j=0} p_j  \frac{\lambda}{\mu_j} ( 2 - e^{-\mu_j t} )  .
\end{eqnarray}

\begin{proof}
The proof follows easily from conditioning on the hyper-exponential state.  
\end{proof}
\end{proposition}

\begin{figure}[ht!]
\hspace{-.5in}~\includegraphics[scale=.2]{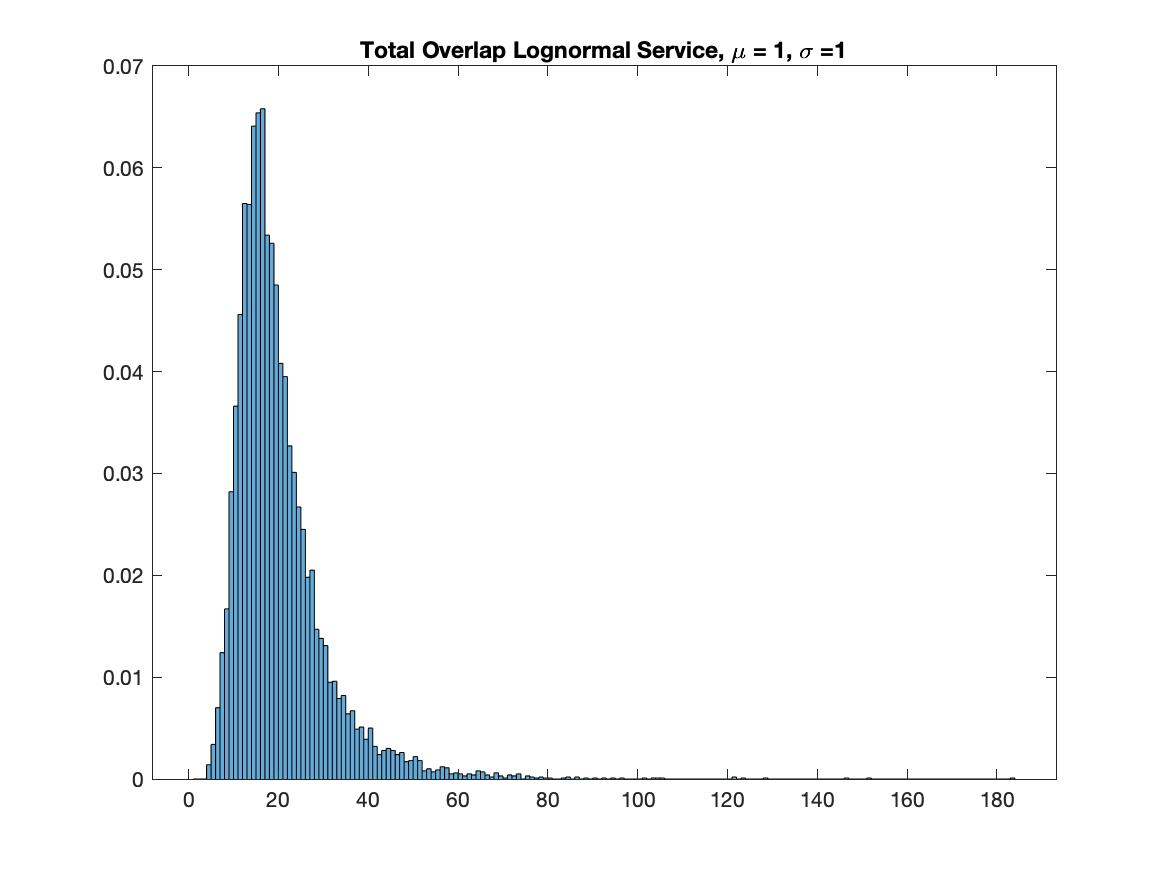}~\hspace{-.25in}~\includegraphics[scale=.2]{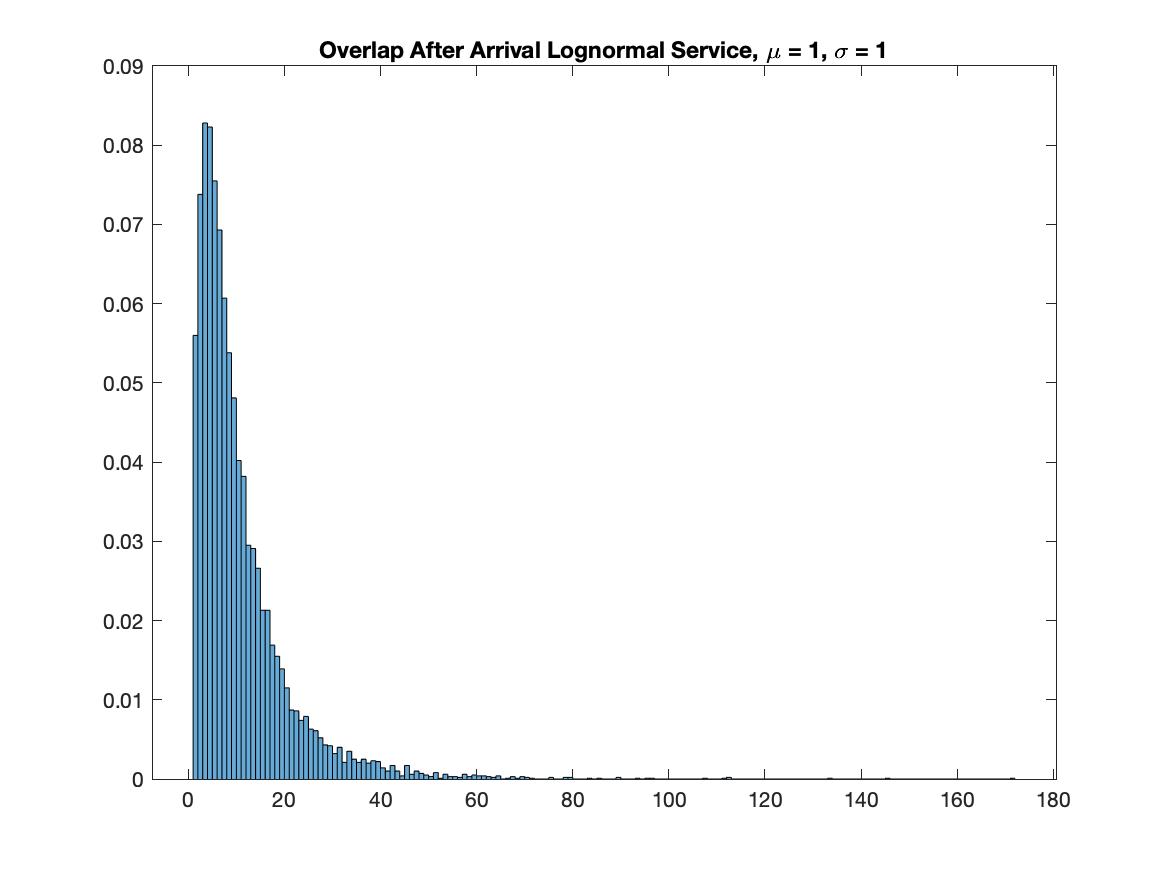}
 \captionsetup{justification=centering}
 \caption{Total Overlap Lognormal Service $(\lambda = 10, \mu = 1, \sigma = 1)$ (Left). \\ 
Overlap During Service Lognormal Service $(\lambda = 10, \mu = 1, \sigma = 1)$ (Right). }
 \label{fig:overlap_all_5}
\end{figure}

In Figure \ref{fig:overlap_all_5}, we plot the overlap distribution and the number of overlaps after arrival for lognormal service times where the mean and variance of the service times are equal to 1.  We see that the total overlap distribution is centered around 20 while the mean of the arrival during service is roughly 10.  One important thing to note is that the number of arrivals during service distribution appears to be similar to an exponential with rate $\mu =1$, however, we observe that the tail of the lognormal is larger.  In fact, we notice that there are a few times where the number overlaps is near 200, while this does not get larger than 50 for the exponential case.  This is clearly connected to the fact that the lognormal distribution does not have a moment generating function all positive values.

\begin{figure}[ht!]
\hspace{-.5in}~\includegraphics[scale=.2]{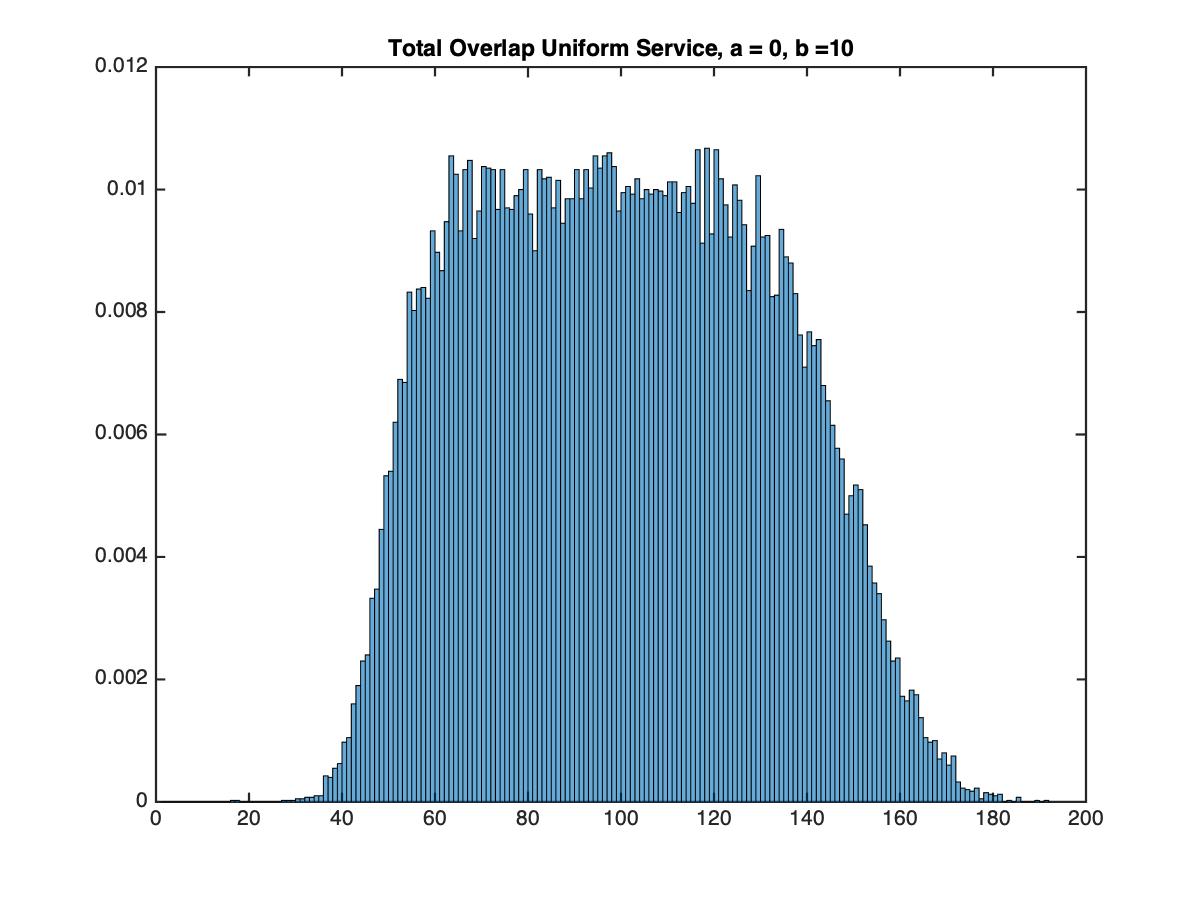}~\hspace{-.25in}~\includegraphics[scale=.2]{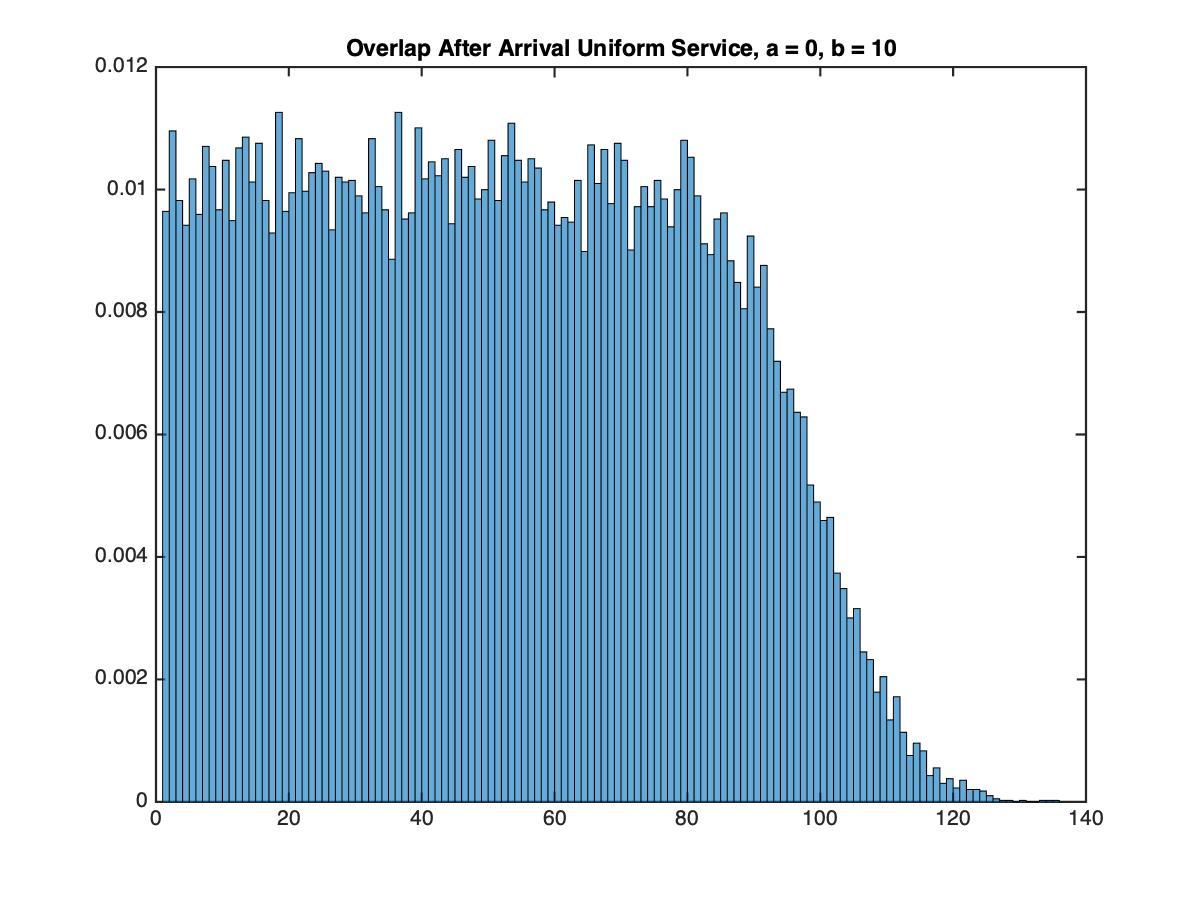}

\hspace{1in}~\includegraphics[scale=.2]{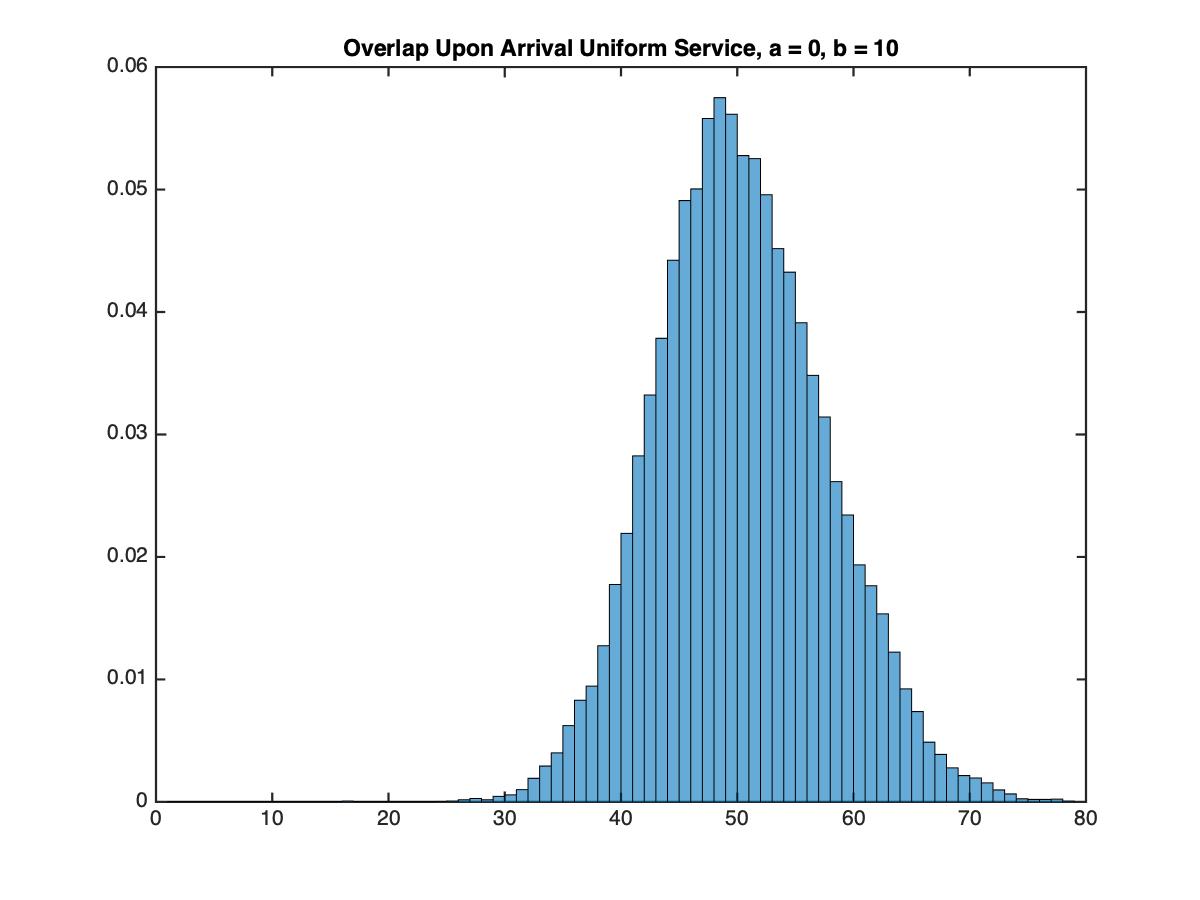}
 \captionsetup{justification=centering}
 \caption{Total Overlap Uniform $(\lambda = 10, a = 0, b = 10)$ (Top Left). \\ 
 Overlap During Service Uniform  $(\lambda = 10, a = 1, b = 10)$ (Top Right). \\
 Overlap Upon Service Uniform  $(\lambda = 10, a = 1, b = 10)$ (Bottom).}
 \label{fig:overlap_all_6}
\end{figure}

In Figure \ref{fig:overlap_all_6}, we plot the overlap distribution and the number of overlaps after arrival for uniform service times in the interval [0,10]. We see that the total overlap distribution is centered around 100 while the mean of the arrival during service is roughly 50.  What is fascinating is that the distribution of overlap is clustered around the mean in a pretty uniform manner.  Moreover, for the overlaps during service, the distribution is quite uniform until it hits 100, which is the boundary of the uniform times the arrival rate, then the distribution begins to decay exponentially, which is similar to the decay obtained by the probability of a Poisson random variable exceeding its mean.  Finally, on the bottom of Figure \ref{fig:overlap_all_6} we plot the overlap upon arrival, which is Poisson with mean roughly equal to 50.  

\section{Residual Overlap Time}

In addition to understanding the distribution of the number of overlapping customers at any point in time, it is also important to understand how long you spend overlapping with those customers as well.  In what follows, we introduce a residual overlap process.  This residual overlap process $O(t,\delta)$ counts the number of customers that you overlap with at least for $\delta$ units of time.  Thus, the process $O(t,\delta)$ has the following representation

\begin{eqnarray}
O(t,\delta) &=& \sum^{N^{(b)}(t + (\mathcal{S} - \delta)^+)}_{j=  N^{(b)}(t) +1} \{ S_j \geq \delta\}  + \sum_{j=1}^{N^{(b)}(t)} \{ ( A_j + S_j - t)^+ \geq \delta \} \cdot \{ \mathcal{S} \geq \delta\} 
\end{eqnarray}

\begin{theorem}
The mean and variance of the residual overlap time is given by the following expressions
\begin{eqnarray}
E[O(t,\delta)] &=&  \lambda E[(\mathcal{S} - \delta)^+ ] \cdot \overline{G}(\delta)  + \lambda \overline{G}(\delta)\int^{t}_{0} \overline{G}(t+\delta-u) du \\
\mathrm{Var}[O(t,\delta)] &=&  \left( \lambda^2 \mathrm{Var} \left[ (\mathcal{S} - \delta)^+  \right] + \lambda E[(\mathcal{S} - \delta)^+ ] \right) \cdot \overline{G}(\delta)^2 + \overline{G}(\delta) \cdot G(\delta) \cdot \lambda E[(\mathcal{S} - \delta)^+ ]  \nonumber \\ \nonumber 
&+&   \left( \lambda \int^{t}_{0} \overline{G}(t+\delta-u) du\right) \cdot \left( \overline{G}(\delta) \cdot G(\delta) + \overline{G}(\delta)^2 \right) \\ 
&+& \overline{G}(\delta) \cdot G(\delta) \cdot \left( \lambda \int^{t}_{0} \overline{G}(t+\delta-u) du \right)^2   .
\end{eqnarray}
\begin{proof}
For the mean we have the following derivation
\begin{eqnarray*}
E[O(t,\delta)] &=&  \mathbb{E} \left[ \sum^{N^{(b)}(t + (\mathcal{S} - \delta)^+)}_{j=  N^{(b)}(t) +1} \{ S_j \geq \delta\} \right]  +  \mathbb{E} \left[  \sum_{j=1}^{N^{(b)}(t)} \{ ( A_j + S_j - t)^+ \geq \delta \} \cdot \{ \mathcal{S} \geq \delta\} \right]  \\
&=&  \mathbb{E} \left[ N^{(b)}(t + (\mathcal{S} - \delta)^+) - N^{(b)}(t) \right] \cdot  \mathbb{E} \left[\{ S_j \geq \delta\} \right]  \\
&+&  \mathbb{E} \left[  \sum_{j=1}^{N^{(b)}(t)} \{ ( A_j + S_j - t)^+ \geq \delta \} \right] \cdot \mathbb{E} \left[\{ \mathcal{S} \geq \delta\} \right]  \\
&=&  \lambda E[(\mathcal{S} - \delta)^+ ] \cdot \overline{G}(\delta) + \lambda  \overline{G}(\delta)  \int^{t}_{0} \overline{G}(t+\delta-u) du.
\end{eqnarray*}

Finally, for the variance we have
\begin{eqnarray*}
\mathrm{Var}[O(t,\delta)]  &=&  \mathrm{Var} \left[ \sum^{N^{(b)}(t + (\mathcal{S} - \delta)^+)}_{j=  N^{(b)}(t) +1} \{ S_j \geq \delta\} \right]  +   \mathrm{Var}  \left[  \sum_{j=1}^{N^{(b)}(t)} \{ ( A_j + S_j - t)^+ \geq \delta \} \cdot \{ \mathcal{S} \geq \delta\} \right]  \\
&=&   \mathrm{Var} \left[ \sum^{N^{(b)}((\mathcal{S} - \delta)^+)}_{j=1} \{ S_j \geq \delta\} \right]  +   \mathrm{Var}  \left[  Q^{\infty}(t,\delta) \cdot \{ \mathcal{S} \geq \delta\} \right]  \\
&=&   \mathrm{Var} \left[ \sum^{N^{(b)}((\mathcal{S} - \delta)^+)}_{j=1} \{ S_j \geq \delta\} \right]  +   \mathrm{Var}  \left[  Q^{\infty}(t,\delta) \right] \cdot  \mathrm{Var} \left[ \{ \mathcal{S} \geq \delta\} \right] \\
&+&   \mathrm{Var}  \left[  Q^{\infty}(t,\delta) \right] \cdot E\left[ \{ \mathcal{S} \geq \delta\} \right]^2 +  \mathrm{Var} \left[ \{ \mathcal{S} \geq \delta\} \right] \cdot E \left[  Q^{\infty}(t,\delta) \right]^2   \\
&=&  \left( \lambda^2 \mathrm{Var} \left[ (\mathcal{S} - \delta)^+  \right] + \lambda E[(\mathcal{S} - \delta)^+ ] \right) \cdot \overline{G}(\delta)^2 + \overline{G}(\delta) \cdot G(\delta) \cdot \lambda E[(\mathcal{S} - \delta)^+ ] \\
&+&   \mathrm{Var}  \left[  Q^{\infty}(t,\delta) \right] \cdot \overline{G}(\delta) \cdot G(\delta)  \\
&+&   \mathrm{Var}  \left[  Q^{\infty}(t,\delta) \right] \cdot \overline{G}(\delta)^2 +  \overline{G}(\delta) \cdot G(\delta) \cdot E \left[  Q^{\infty}(t,\delta) \right]^2   \\
&=&  \left( \lambda^2 \mathrm{Var} \left[ (\mathcal{S} - \delta)^+  \right] + \lambda E[(\mathcal{S} - \delta)^+ ] \right) \cdot \overline{G}(\delta)^2 + \overline{G}(\delta) \cdot G(\delta) \cdot \lambda E[(\mathcal{S} - \delta)^+ ] \\
&+&   \mathrm{Var}  \left[  Q^{\infty}(t,\delta) \right] \cdot \left( \overline{G}(\delta) \cdot G(\delta) + \overline{G}(\delta)^2 \right) + \overline{G}(\delta) \cdot G(\delta) \cdot E \left[  Q^{\infty}(t,\delta) \right]^2  .
\end{eqnarray*}
This completes the proof.
\end{proof}
\end{theorem}

Using the mean and variance of the residual overlap process, one can construct an approximate prediction interval for how many people are expected to overlap with a customer that arrives at time $t$. This analysis is useful since it can inform managers of how to "restrict" the arrival rate to achieve a particular mean overlap or type of distribution of overlap that is sufficient for risk tolerance for employees and patrons.

\begin{proposition}
Let $\mathcal{S}$ be a continuous and non-negative random variable with cdf given by $G(x)$ and define $N((\mathcal{S} - \delta)^+)$ to be the the number of arrivals during service given that service, then its distribution is given by the following expression
\begin{eqnarray*}
\mathbb{P} \left( N((\mathcal{S} - \delta)^+) = k \right) &=& \frac{\lambda^k}{k!} e^{\lambda \delta} \mathbb{E} \left[ (\mathcal{S} - \delta)^k e^{-\lambda \mathcal{S}} \bigg | \mathcal{S} > \delta \right] \overline{G}(\delta) \\
\mathbb{P} \left( N((\mathcal{S} - \delta)^+) = 0 \right)&=& G(\delta) + e^{\lambda \delta} \mathbb{E}\left[ e^{-\lambda \mathcal{S}} \bigg | \mathcal{S} > \delta \right] \overline{G}(\delta) .
\end{eqnarray*}
\begin{proof}
\begin{eqnarray*}
\mathbb{P} \left( N((\mathcal{S} - \delta)^+)  = k  \right)  &=&  \int^{\infty}_{\delta } e^{-\lambda (x-\delta)} \frac{(\lambda (x-\delta))^k}{k!} dG(x) \\
&=& \frac{\lambda^k}{k!} e^{\lambda \delta} \int^{\infty}_{\delta } (x-\delta)^k e^{-\lambda x}  dG(x) \\
&=& \frac{\lambda^k}{k!} e^{\lambda \delta} \mathbb{E} \left[ (\mathcal{S} - \delta)^k e^{-\lambda \mathcal{S}} \bigg | \mathcal{S} > \delta \right] \overline{G}(\delta) .
\end{eqnarray*}
The term for zero follows easily by computing what is left over so the sum totals one.  
\end{proof}
\end{proposition}

\begin{corollary}
Let $\mathcal{S}$ be a Exp$(\mu)$ random variable, then 
\begin{eqnarray}
\mathbb{P} \left( N((\mathcal{S}- \delta)^+) = k \right) =  \left( \frac{\mu }{\lambda + \mu} \right) \left( \frac{ \lambda}{\lambda+\mu} \right)^k e^{-\mu \delta} , \quad k \geq 1  \\
\mathbb{P} \left( N((\mathcal{S}- \delta)^+)  = 0 \right) = 1 -  \left( \frac{\lambda }{\lambda + \mu} \right)  e^{-\mu \delta}  . 
\end{eqnarray}
\begin{proof}
This follows by directly computing the necessary conditional expectation. 
\begin{eqnarray*}
\mathbb{P} \left( N((\mathcal{S} - \delta)^+)  = k  \right)  &=&  \int^{\infty}_{\delta } e^{-\lambda (x-\delta)} \frac{(\lambda (x-\delta))^k}{k!} \mu e^{-\mu x} dx \\
&=& \frac{\mu \lambda^k}{k!} e^{\lambda \delta} \int^{\infty}_{\delta } (x-\delta)^k e^{-(\lambda+\mu) x} dx \\
&=& \frac{\mu \lambda^k}{k!} e^{-\mu \delta} \int^{\infty}_{0} y^k e^{-(\lambda+\mu) y} dx \\
&=& \frac{\mu \lambda^k}{k!} e^{-\mu \delta} \frac{k!}{(\lambda + \mu)^{k+1}} \\
&=& \left( \frac{\mu }{\lambda + \mu} \right) \left( \frac{ \lambda}{\lambda+\mu} \right)^k e^{-\mu \delta} .
\end{eqnarray*}
The remaining term when $k=0$ follows from summing all of the other terms and subtracting it from 1.  This completes the proof.  
\end{proof}
\end{corollary}

\begin{proposition}
Let $\mathcal{S}$, be a Exp$(\mu)$ random variable, $S_j$ i.i.d random variables with cdf $G(x)$ and define $Z(t,\delta)$ to be the the number of arrivals during service that overlap at least $\delta$ units of time i.e.
\begin{eqnarray}
Z(t,\delta) = \sum^{N^{(b)}(t + (\mathcal{S} - \delta)^+)}_{j=  N^{(b)}(t) +1} \{ S_j \geq \delta\} .
\end{eqnarray}
Then, we have the following expression for the distribution of $Z(t,\delta)$
\begin{eqnarray*}
\mathbb{P} \left( Z(t,\delta)  = k  \right)  &=&   ( 1 - \rho)\rho^{k}  e^{-\mu \delta} \overline{G}(\delta)^{k}  \frac{1}{(1- \rho G(\delta))^{k+1}} .
\end{eqnarray*}
\begin{proof}
\begin{eqnarray*}
\mathbb{P} \left( Z(t,\delta)  = k  \right)  &=& \sum^\infty_{\ell=0} \mathbb{P} \left( Z(t,\delta)  = k  \bigg | N^{(b)}( (\mathcal{S} - \delta)^+) = \ell \right) \cdot  \mathbb{P} \left(  N^{(b)}( (\mathcal{S} - \delta)^+) = \ell \right) \\
&=& \sum^\infty_{\ell=k} { \ell \choose k } \overline{G}(\delta)^k G(\delta)^{\ell - k} \cdot  ( 1 - \rho) \rho^{\ell}  e^{-\mu \delta} \\
&=&  \left( \frac{\overline{G}(\delta) }{G(\delta)} \right)^{k} ( 1 - \rho) e^{-\mu \delta}  \sum^\infty_{\ell=k} { \ell \choose k }  G(\delta)^{\ell } \cdot   \rho^{\ell}  \\
&=&  \left( \frac{\overline{G}(\delta) }{G(\delta)} \right)^{k} ( 1 - \rho) e^{-\mu \delta}  \sum^\infty_{\ell=k} { \ell \choose k }  (\rho G(\delta))^{\ell }  \\
&=&  \left( \frac{\overline{G}(\delta) }{G(\delta)} \right)^{k} ( 1 - \rho) e^{-\mu \delta}  \sum^\infty_{j=0} { j+k \choose k }  (\rho G(\delta))^{j+k }  \\
&=&  \left( \overline{G}(\delta) \rho \right)^{k} ( 1 - \rho) e^{-\mu \delta}  \sum^\infty_{j=0} { j+k \choose k }  (\rho G(\delta))^{j}  \\
&=&  \left( \overline{G}(\delta) \rho \right)^{k} ( 1 - \rho) e^{-\mu \delta}  \frac{1}{(1- \rho G(\delta))^{k+1}} .
\end{eqnarray*}
\end{proof}
\end{proposition}

\begin{figure}[ht!]
\hspace{-.5in}~\includegraphics[scale=.2]{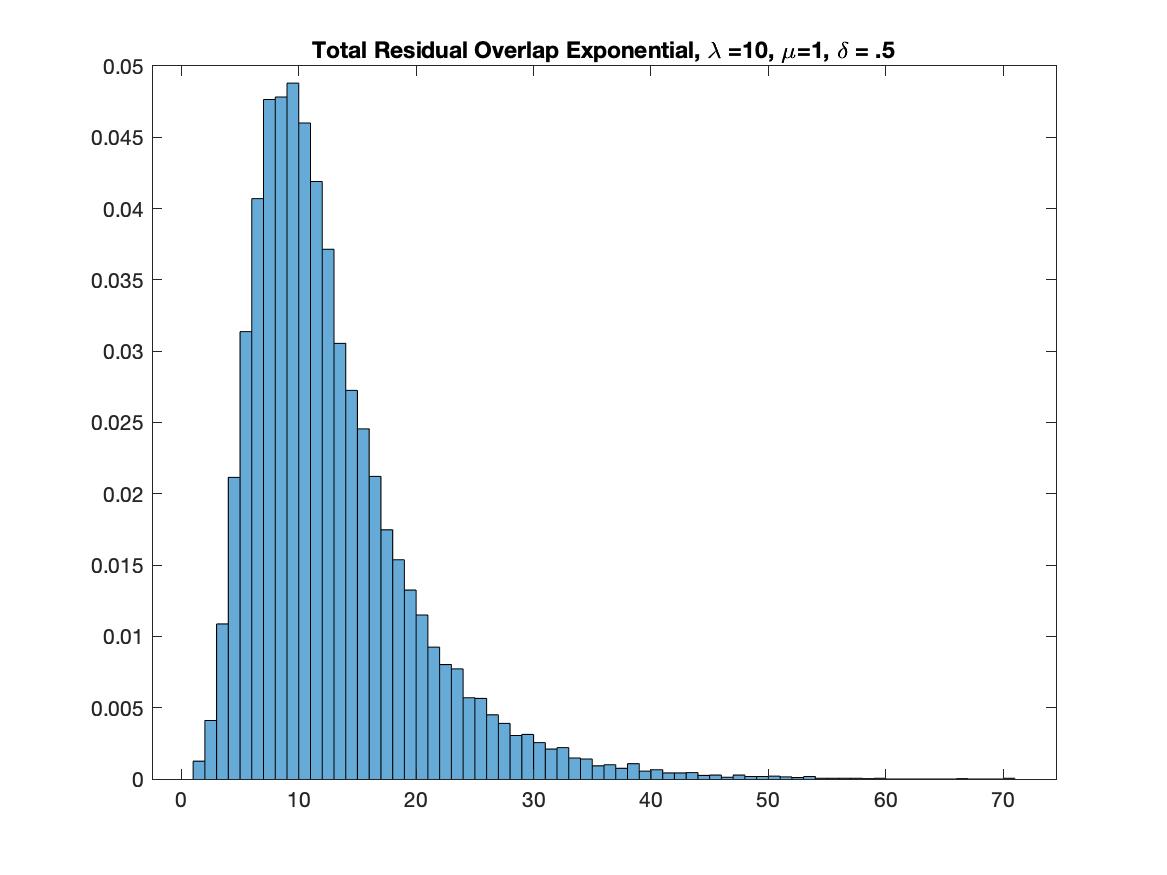}~\hspace{-.25in}~\includegraphics[scale=.2]{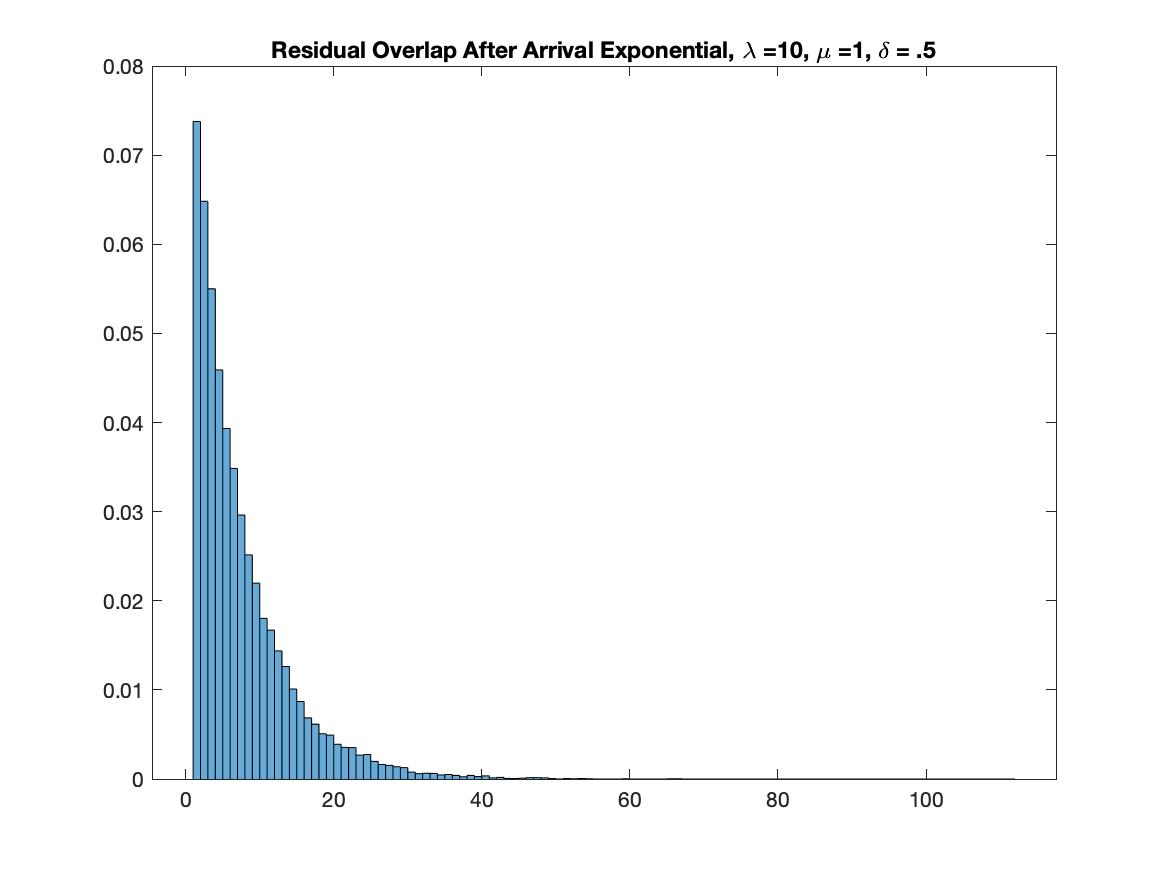}

\hspace{1in}~\includegraphics[scale=.2]{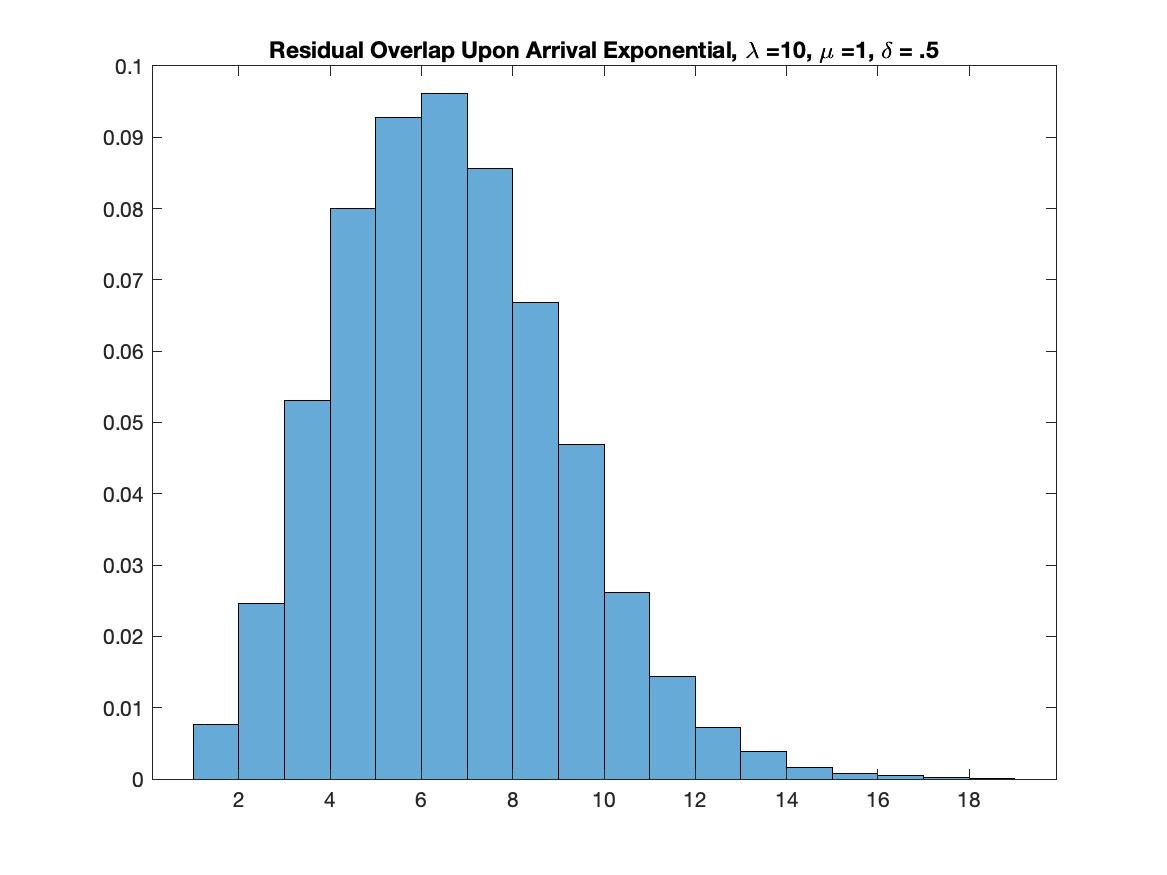}
 \captionsetup{justification=centering}
 \caption{Total Residual Overlap (Exponential), $(\lambda = 10, \mu = 1, \delta = .5)$ (Top Left). \\ 
 Residual Overlap After Arrival (Exponential), $(\lambda = 10, \mu = 1, \delta =.5)$ (Top Right). \\
 Residual Upon Arrival (Exponential), $(\lambda = 10, \mu = 1, \delta =.5)$ (Bottom). }
 \label{fig:overlap_all_7}
\end{figure}

In Figure \ref{fig:overlap_all_7}, we plot the residual distribution for the total number of overlaps, the number of overlaps upon arrival and the number overlaps during service.  We see that the distribution is much different than that of Figure \ref{fig:overlap_all_8}.  We see much fewer overlaps since we are requiring at least an overlap of .5 time units, which is decently large given the service time has a mean equal to one.  


\section{Conclusion and Future Work} \label{conclusion}

In this paper, we consider the overlap times for customers in an infinite server queue.  The infinite server model is appropriate in retail settings where the time a customer waits is small relative to their shopping experience.  We derive the steady state distribution for the overlap time of customers that are $k$ arrivals apart.  We also compute explicitly the distribution of the number of customers that a randomly arriving customer will overlap with as the sum of a Poisson random variable and Poisson with a random arrival rate.  Finally, we compute an expression for the number of customers that a random arrival will overlap at least $\delta$ time units.  We compute the mean and variance for this quantity and are able to provide a prediction interval for a randomly arriving customer at any time.  Our analysis has implications for understanding the interaction time between customers in a pandemic setting and sheds light on the interactions of customers.   

Despite our analysis, there are many avenues for additional research.  First, we would like to complete our analysis by analyzing the $G/G/\infty$ queue in explicit detail.  Issues like dependent arrivals or service times like in \citet{pang2012impact, pang2013two, daw2018exact, daw2018queues, daw2020co} would be interesting to explore as well.  Second, we would like to extend our analysis to more complicated queueing systems like the Erlang-A, see for example \citet{daw2019new, massey2018dynamic, hampshire2020beyond} where customers can abandon the system.  Abandoning customers clearly reduces the number of overlapping customers, but by how much?  It would also be great to extend our analysis to queueing systems with batch arrivals.  In this case, the scaled Poisson decomposition of \citet{daw2019distributions,daw2020non} might be helpful in replicating our analysis in the batch setting.  Finally, we would like to extend our analysis to spatial point processes and think about the overlap in terms of not only time, number, but spatial distance as well.  

We are also interested in knowing the overlap distribution for more than two customers.  For example, if one considers the overlap time of three customers $(n,n+j,n+k)$ where $1 \leq j < k $, then one obtains the following overlap time for the $n^{th}$, $(n+j)^{th}$, and the $(n+k)^{th}$ customers as

\begin{eqnarray*}
O_{n,n+j,n+k} &=& \left( \min( D_n , D_{n+j}, D_{n+k} ) - A_{n+k} \right)^+ \\
&=& \left( \min( A_n + S_n , A_{n+j} + S_{n+j},  A_{n+k} + S_{n+k} ) - A_{n+k} \right)^+ \\
&=& (S_n + A_{n} - A_{n+k} )^+ \boldsymbol{\wedge} (S_{n+j} + A_{n+j} - A_{n+k} )^+ \boldsymbol{\wedge} S_{n+k}  \\
&=& (S_n + A_{n} - A_{n+k} )^+ \boldsymbol{\wedge} (S_{n+j} + A_n + (A_{n+j} - A_n) - A_{n+k} )^+ \boldsymbol{\wedge} S_{n+k}  .
\end{eqnarray*}
Unlike the two customer situation, it is clear here that the random variables are not independent anymore. This presents a new issue that must be resolved in future work.

\bibliographystyle{plainnat}
\bibliography{infinite_overlap}
\end{document}